\documentclass{amsart}

\usepackage{amssymb,amsmath,xcolor}
\usepackage[
  colorlinks=true,
  linkcolor=red!70!black,
  citecolor=blue!75!black,
  urlcolor=blue!75!black,
  pdfusetitle
]{hyperref}

\newtheorem{theorem}{Theorem}[section]
\newtheorem{proposition}[theorem]{Proposition}
\newtheorem{lemma}[theorem]{Lemma}
\newtheorem{corollary}[theorem]{Corollary}
\theoremstyle{remark}
\newtheorem{remark}[theorem]{Remark}
\newtheorem{step}{Step}
\theoremstyle{definition}
\newtheorem{definition}[theorem]{Definition}
\newtheorem{condition}[theorem]{Condition}
\newtheorem{example}[theorem]{Example}

\DeclareMathOperator{\Var}{Var}
\DeclareMathOperator{\Cov}{Cov}
\DeclareMathOperator{\sech}{sech}
\DeclareMathOperator{\artanh}{artanh}
\DeclareMathOperator{\sgn}{sgn}
\DeclareMathOperator{\Real}{Re}
\DeclareMathOperator{\Img}{Im}
\DeclareMathOperator{\dist}{dist}
\newcommand{\E}{\mathbb E}
\newcommand{\R}{\mathbb R}
\newcommand{\Z}{\mathbb Z}
\newcommand{\wt}{\widetilde}
\newcommand{\HS}{\mathrm{HS}}

\def\MR#1{\href{https://mathscinet.ams.org/mathscinet-getitem?mr=#1}{MR#1}}

\numberwithin{equation}{section}

\title[Real zeros of generalized Kac polynomials]{Logarithmic variance of the number of real zeros of generalized Kac polynomials}

\author{Nhan D. V. Nguyen}

\address{Department of Mathematics and Statistics, Binghamton University (SUNY), Binghamton, NY 13902-6000 \and Department of Mathematics and Statistics, Quy Nhon University, Gia Lai, Vietnam}

\email{nnguyen36@binghamton.edu, nguyenduvinhan@qnu.edu.vn}

\keywords{Random polynomials, real zeros, variance, central limit theorem, Kac--Rice formula, Gaussian processes, Wiener chaos}

\subjclass[2020]{30C15, 60F05, 60G15, 26C10}
\date{\today}

\begin{document}
\begin{abstract}
    Let $N_n$ be the number of real zeros of the random polynomial $P_n(x)=\sum_{j=0}^n\xi_jv_jx^j$, where the random coefficients $\xi_j$ are independent, centered, with unit variance and uniformly bounded $(2+\varepsilon_0)$-moments, and the deterministic weights satisfy $|v_j|\asymp(1+j)^\tau$ with $\tau>-1/2$, uniformly in $n$ and $j$. No further regularity of the weights is assumed. We prove that $\Var N_n\asymp\log n$. Consequently, the variance hypothesis in the central limit theorem of O.~Nguyen and Vu holds for all such weights, so $N_n$ satisfies a central limit theorem. After a reduction to Gaussian coefficients on a core region near $\pm1$, the upper bound follows from Jensen's formula and the decorrelation of the two-point zero density at large hyperbolic distance, and the lower bound from the projection of $N_n$ onto the second Wiener chaos, tested on Mellin-modulated vectors. Using the comparison principles of Do, Nguyen, and O'Rourke, we extend the variance bounds and the central limit theorem to $\sum_{j=0}^n(\mu_j+\xi_j)v_jx^j$ with bounded $\mu_j$, under suitable domination conditions on the mean.
\end{abstract}
\maketitle

\section{Introduction}
    How many real solutions does a random algebraic equation have? By a classical result of Kac \cite{Kac43}, a polynomial of degree $n$ with independent standard Gaussian coefficients has on average $[2/\pi+o(1)]\log n$ real zeros, far fewer than its $n$ complex zeros. These real zeros accumulate near $\pm1$. Heuristically, in the hyperbolic variable $t=\artanh x$, the region $1/n\le 1-x\le 1$ corresponds to an interval of length of order $\log n$, on which the normalized process has speed of order one and its correlations decay exponentially in the hyperbolic distance. The zero count therefore behaves like a sum of order $\log n$ weakly dependent contributions, and its variance is also of order $\log n$. In this paper, we determine the order of the variance of the number of real zeros for a broad class of polynomials with independent coefficients and non-identical weights, without any regularity assumption on the weights.

    Consider the random polynomial
\begin{equation} \label{eq:poly}
P_n(x)=\sum_{j=0}^n \xi_jv_jx^j,
\end{equation}
    where $\xi_0, \dots, \xi_n$ are independent, real-valued random variables and $v_0, \dots, v_n$ are deterministic real weights. Let $N_n$ denote the number of real zeros of $P_n$, counted with multiplicity, with the convention that $N_n=0$ on the event $\{P_n\equiv 0\}$. We study the fluctuations of $N_n$ as $n\to\infty$.
    
\subsection{Background}
    When $v_j\equiv 1$, $P_n$ is the classical Kac polynomial. The first estimates for the number of its real zeros are due to Bloch and P\'olya \cite{BP32}, who bounded the average number of real zeros by $O(\sqrt n)$ for coefficients uniformly distributed on $\{-1,0,1\}$. Littlewood and Offord \cite{LO43, LO45, LO48} later improved this to polylogarithmic order. For Gaussian coefficients, Kac \cite{Kac43} derived an integral formula for the expected density of real zeros, an early form of the Kac--Rice method \cite{AW09, EK95}, and proved that
\[
\E N_n=\left[\frac{2}{\pi}+o(1)\right]\log n.
\]
    This asymptotic formula persists for many coefficient distributions \cite{EO56, Kac49, St69, IM71P1}. Universality results \cite{DHNV15, NNV16, NV22, TV15} explain this stability by comparing local zero statistics for different coefficient laws under suitable moment and non-degeneracy assumptions.

\begin{condition}\label{cond:xi}
    The random variables $\xi_j$, $j\ge0$, are independent and real-valued, with 
\[
\E\xi_j=0,\quad \Var\xi_j=1, \quad \text{and}\quad \sup_{j\ge 0}\E|\xi_j|^{2+\varepsilon_0}\le C_0,
\]
    for some constants $C_0,\varepsilon_0>0$.
\end{condition}

    Maslova \cite{M74,M74I} obtained the first results on the fluctuations of $N_n$. For i.i.d.\ coefficients satisfying Condition \ref{cond:xi} and $\mathbb{P}(\xi_j=0)=0$, she proved that
\[
\Var N_n = \left[\frac{4}{\pi}\left(1-\frac{2}{\pi}\right)+o(1)\right]\log n
\]
    and established the central limit theorem (CLT)
\[
\frac{N_n-\E N_n}{\sqrt{\Var N_n}}
\xrightarrow{d}\mathcal N(0,1).
\]
    Her argument approximates $N_n$ by a sum of weakly dependent local zero counts over carefully chosen intervals. Similar results have since been obtained for other classes of random polynomials and random functions, in both Gaussian and non-Gaussian settings. These include random trigonometric polynomials \cite{ADL16, AL13, BCP19, DNN22, Gas23}, Kostlan, Weyl, and elliptic polynomials \cite{AN26, ANW26, AL21, BD97, Dal15, DV20, N24, SM08}, random orthogonal polynomials \cite{DNNP24, LP21, LP22}, and stationary Gaussian processes \cite{AL21b, ABF23}. 

    This paper focuses on generalized Kac polynomials, in which the weights $v_j$ vary polynomially in $j$. This class contains the Kac polynomials, the hyperbolic polynomials, and finite linear combinations of their derivatives \cite{DN25,NV21}.

\begin{condition}\label{cond:main}
    There exist constants $0<C_1\le C_2$ and $\tau>-1/2$ such that, for every $n\ge 1$ and $0\le j\le n$, the weights $v_j=v_{j,n}$ are real and satisfy
\begin{equation}\label{eq:weights}
C_1(1+j)^\tau\le|v_j|\le C_2(1+j)^\tau.
\end{equation} 
\end{condition}

    Condition \ref{cond:main} allows the $v_j$ to have arbitrary signs, and allows $|v_j|$ to oscillate between $C_1(1+j)^\tau$ and $C_2(1+j)^\tau$, for instance along the dyadic scales of $j$. For $\tau=0$, it holds for every choice of weights with $C_1\le|v_j|\le C_2$.

    The expected number of real zeros is well understood in this setting. Under Conditions \ref{cond:xi} and \ref{cond:main}, Do, Nguyen, and Vu \cite{DNV18} proved that
\[
c \log n \le \E N_n \le C \log n.
\]
    If, in addition, $|v_j|\sim A(1+j)^\tau$ for some constant $A>0$, they obtained
\[
\E N_n = \left[\frac{1+\sqrt{2\tau+1}}{\pi}+o(1)\right] \log n.
\]

    The fluctuations are more delicate because varying weights make it harder to decompose $N_n$ into weakly dependent local contributions. Nguyen and Vu \cite{NV21} combined universality arguments with a refinement of Maslova's method to prove a CLT for generalized Kac polynomials under a logarithmic lower bound on the variance.

\begin{theorem}[{\cite[Theorem 1.2]{NV21}}]\label{thm:NV}
    Let $P_n$ be defined by \eqref{eq:poly}, where Conditions \ref{cond:xi} and \ref{cond:main} hold. If
\begin{equation}\label{eq:lowbound-ass}
\Var N_n\ge c\log n\quad\text{for some }c>0\text{ and all sufficiently large }n,
\end{equation}
    then
\[
\frac{N_n-\E N_n}{\sqrt{\Var N_n}}
\xrightarrow{d}\mathcal N(0,1).
\]
\end{theorem}

    Nguyen and Vu verified \eqref{eq:lowbound-ass} under the additional assumption that, for some $\varepsilon>0$,
\begin{equation} \label{eq:NV-cond}
    \frac{|v_j|}{|v_n|}-1=O\left(e^{-(\log\log n)^{1+\varepsilon}}\right),\quad  n-ne^{-\log^{1/5}n}\le j \le n-e^{\log^{1/5}n}.
\end{equation}
    This holds for Kac polynomials, hyperbolic polynomials, and their derivatives.

    Under the alternative regularity assumption $|v_j|\sim A(1+j)^\tau$, Do and the author \cite{DN25} obtained the sharp asymptotic formula
\[
\Var N_n=\left[2\beta_\tau+\frac2\pi\left(1-\frac2\pi\right)+o(1)\right]\log n,
\]
    where $\beta_\tau$ is an explicit constant depending only on $\tau$, with $\beta_0= \pi^{-1}(1-2/\pi)$.

\subsection{Main results}
    By Theorem \ref{thm:NV}, the central limit theorem for $N_n$ follows from \eqref{eq:lowbound-ass}. This bound was previously known under the additional condition \eqref{eq:NV-cond} \cite{NV21} and under the asymptotic assumption $|v_j|\sim A(1+j)^\tau$ \cite{DN25}. For weights satisfying only \eqref{eq:weights}, the normalized covariance kernels need not converge to the limiting kernels used in \cite{NV21, DN25}, and the arguments there do not apply directly. To the best of our knowledge, matching upper and lower bounds of order $\log n$ for $\Var N_n$ were not previously available for this class of weights.

\begin{theorem}[Logarithmic variance] \label{thm:main}
    Let $P_n$ be defined by \eqref{eq:poly}, where Conditions \ref{cond:xi} and \ref{cond:main} hold. Then there exist constants $C, c>0$, depending only on $\tau,\varepsilon_0,C_0,C_1,C_2$, such that for all sufficiently large $n$,
\[
c\log n\le\Var N_n\le C\log n.
\]
\end{theorem}

\begin{remark} 
    Theorem \ref{thm:main} determines the order of $\Var N_n$ and does not assert that $\Var N_n/\log n$ converges. The point is that no regularity of the weights beyond \eqref{eq:weights} is required. The proof replaces an explicit limiting covariance kernel by a localization argument for the upper bound and by Mellin-modulated test vectors for the lower bound.
\end{remark}

    The lower bound in Theorem \ref{thm:main} verifies the hypothesis \eqref{eq:lowbound-ass} of Theorem \ref{thm:NV} and therefore yields the following central limit theorem. 
    
\begin{corollary}[CLT]\label{cor:clt}
    Under the hypotheses of Theorem \ref{thm:main}, as $n\to \infty$,
\[
\frac{N_n-\E N_n}{\sqrt{\Var N_n}} \xrightarrow{d} \mathcal N(0,1).
\]
\end{corollary}

\subsection{Non-centered polynomials} \label{subsec:noncentered}
    We now extend Theorem \ref{thm:main} and Corollary \ref{cor:clt} to non-centered generalized Kac polynomials, provided that the mean function satisfies suitable domination conditions. This extension follows from Theorem \ref{thm:main} and the comparison theorems of \cite{DNOR26}.

    For the classical Kac ensemble with i.i.d.\ non-centered coefficients, Ibragimov and Maslova \cite{IM71II} showed that the leading asymptotics of the expected number of real zeros is one-half of that in the centered case. Do \cite{Do21} extended this result to generalized Kac polynomials by comparing the deterministic mean with the centered random part. This comparison identifies regions where the deterministic mean suppresses real zeros and regions where the centered random fluctuation determines their number. More recently, Do, Nguyen, and O'Rourke \cite{DNOR26} extended these comparison principles to second-order fluctuations.

    Let
\begin{equation}\label{eq:Qn}
Q_n(x)=\sum_{j=0}^n(\mu_j+\xi_j)v_jx^j=M_n(x)+P_n(x),\qquad M_n(x):=\sum_{j=0}^n\mu_jv_jx^j,
\end{equation}
    where the $\xi_j$ satisfy Condition \ref{cond:xi}, the $v_j$ satisfy Condition \ref{cond:main}, and the means satisfy $|\mu_j|\le C_\mu$ for all $n\ge 1$, $0\le j\le n$, and some constant $C_\mu>0$. 
    
    We first recall the relevant definitions from \cite{DNOR26}. Let $F^*(x)=x^nF(1/x)$ for a polynomial $F$ of degree at most $n$.
    
\begin{definition}[Domination]
    Let $J\subset\R$ and $\Phi:(0,\infty)\to [0,\infty)$. We say that $M_n$ \textit{dominates} $P_n$ on $J$ with factor function $\Phi$ if
\begin{align*}
|M_n(x)|&\ge\Phi\bigl(1-|x|+1/n\bigr)\sqrt{\Var P_n(x)}, && x\in J\cap[-1,1],\\
|M_n^*(1/x)|&\ge\Phi\bigl(1-|1/x|+1/n\bigr)\sqrt{\Var P_n^*(1/x)}, && x\in J\setminus[-1,1].
\end{align*}
    We say that $M_n$ dominates $P_n$ on $J$ \textit{up to order $2$} if these inequalities hold with $M_n^{(i)}$, $P_n^{(i)}$, $(M_n^*)^{(i)}$, $(P_n^*)^{(i)}$ in place of $M_n$, $P_n$, $M_n^*$, $P_n^*$ for $i=0,1,2$. If the inequalities are reversed, $M_n$ is \textit{dominated by} $P_n$ on $J$ (up to order $2$) with factor function $\Phi$.
\end{definition}

    Roughly speaking, when $M_n$ dominates $P_n$, the deterministic mean is large compared to the typical size of the random fluctuations, and similarly for their derivatives up to the prescribed order.

\begin{definition}[Enlargement]
    An interval $J\supset S$ is an \textit{enlargement} of $S$ if each finite endpoint $x$ of $S$ is extended outward by $c(|1-|x||+1/n)$ for some constant $c>0$, unless $|1-|x||$ is bounded below by a positive constant independent of $n$.
\end{definition}

    For $S\subset\R$, let $N_{Q_n}(S)$ and $N_{P_n}(S)$ be the numbers of real zeros in $S$ of $Q_n$ and $P_n$, counted with multiplicity. As for $P_n$, we set $N_{Q_n}(S)=0$ on the event $\{Q_n\equiv0\}$. We write $N_n(S)=N_{P_n}(S)$ and $N_{Q_n}=N_{Q_n}(\R)$.  The two comparison theorems of \cite{DNOR26} are as follows.

\begin{theorem}[{\cite[Theorems~1.10 and~1.12]{DNOR26}}]\label{thm:DNOR} 
    Let $Q_n$ be defined by \eqref{eq:Qn}, with $\xi_j$ satisfying Condition \ref{cond:xi}, $v_j$ satisfying Condition \ref{cond:main}, and $|\mu_j|\le C_\mu$ for some constant $C_\mu>0$. Let $S$ be an interval (possibly depending on $n$), let $J$ be an enlargement of $S$, and fix $d\in (0,1)$.
\begin{enumerate}
    \item If $M_n$ dominates $P_n$ on $J$ with factor function $y\mapsto \alpha |\log y|^{1/2}$ for $\alpha\ge 2$, then
    \[
    \Var N_{Q_n}(S)=O(\log^{d}n).
    \]
    \item Fix $\beta , \gamma, \epsilon>0$. If $M_n$ is dominated by $P_n$ on $J$ up to order $2$ with factor function $y\mapsto \beta y^{\gamma}$, then
    \[
    \Var N_{Q_n}(S)=\Var N_{P_n}(S)+O(\log^{d} n).
    \]
    If, in addition, $\Var N_{P_n}(S)\ge\epsilon\log n$, then
    \[
    \frac{N_{Q_n}(S)-\E N_{Q_n}(S)}{\sqrt{\Var N_{Q_n}(S)}}\xrightarrow{d}\mathcal N(0,1).
    \]
\end{enumerate}
\end{theorem}

    Combined with Theorem \ref{thm:main}, this yields the following.

\begin{theorem}[Non-centered polynomials]\label{thm:noncentered}
    Let $Q_n$ be defined by \eqref{eq:Qn}, with $\xi_j$ satisfying Condition \ref{cond:xi}, $v_j$ satisfying Condition \ref{cond:main}, and $|\mu_j|\le C_\mu$ for some constant $C_\mu>0$. Fix $\alpha\ge 2$, $\beta,\gamma>0$, and $x_0\in (0,1)$. Let $\Phi_{\alpha}(y)=\alpha |\log y|^{1/2}$ and $\Phi_{\beta,\gamma}(y)=\beta y^{\gamma}$. Suppose that one of the following conditions holds for all sufficiently large $n$.
\begin{enumerate}
\item\label{nc:dominated} 
    $M_n$ is dominated by $P_n$ on $\R$ up to order $2$ with factor function $\Phi_{\beta,\gamma}$.
\item\label{nc:mixed} 
    On one of the intervals $[x_0,1/x_0]$ and $[-1/x_0,-x_0]$, $M_n$ dominates $P_n$ with factor function $\Phi_\alpha$, and on the other it is dominated by $P_n$ up to order $2$ with factor function $\Phi_{\beta,\gamma}$.
\end{enumerate}
    Then there exist constants $C, c>0$, depending only on $\tau$, $\varepsilon_0$, $C_0$, $C_1$, $C_2$, $C_\mu$, $\alpha$, $\beta$, $\gamma$, $x_0$, such that for all sufficiently large $n$,
\[
c\log n\le\Var N_{Q_n}\le C\log n,
\]
    and, as $n\to \infty$, 
\[
\frac{N_{Q_n}-\E N_{Q_n}}{\sqrt{\Var N_{Q_n}}}\ \xrightarrow{d}\ \mathcal N(0,1).
\]
\end{theorem}

    Case \ref{nc:dominated} covers means that are small compared with the fluctuations on all of $\R$. In case \ref{nc:mixed}, the mean suppresses the zeros near one of $\pm1$, while the centered fluctuations determine those near the other. The following two examples with $\tau=0$ illustrate both regimes. 

\begin{example}\label{ex:periodic}
    Let $v_j=1$ if $j\equiv0,1\pmod 4$ and $v_j=-1$ if $j\equiv2,3\pmod 4$, and let $\mu_j=\mu$ for fixed $\mu\in \R$. Then $M_n$ is dominated by $P_n$ on $\R$ up to order $2$ with factor function $\Phi_{\beta,1/2}$ for some $\beta>0$. Thus, case \ref{nc:dominated} of Theorem \ref{thm:noncentered} applies.
\end{example}

\begin{example}\label{ex:alternating}
    Let $v_j=(-1)^j$ and $\mu_j=(-1)^j\mu$ for fixed $\mu\ne 0$. Then for any $\alpha>0$, there exist $x_0\in (0,1)$ and $\beta>0$ such that $M_n$ is dominated by $P_n$ on $[-1/x_0, -x_0]$ up to order $2$ with factor function $\Phi_{\beta,1/2}$ and, for all sufficiently large $n$, $M_n$ dominates $P_n$ on $[x_0, 1/x_0]$ with factor function $\Phi_\alpha$. Thus, case \ref{nc:mixed} of Theorem \ref{thm:noncentered} applies.
\end{example}

    Both examples are verified in Section \ref{sec:noncentered}.
    
\subsection{Strategy of the proof}
    The proof reduces to Gaussian coefficients and then studies the normalized process $P_n(x)/\sqrt{\Var P_n(x)}$ in the hyperbolic variable $t=\artanh x$. Since $1-\tanh t\asymp e^{-2t}$, the core region near $1$ corresponds to an interval $J_n$ of length of order $\log n$. Heuristically, the process decorrelates at distance $O(1)$ in $t$, so $N_n$ behaves like a sum of order $\log n$ weakly dependent contributions. The difficulty is to make this rigorous without stationarity and without a limiting covariance kernel.  In the variable $t$, \eqref{eq:weights} alone ensures that the normalized coefficient curve has speed of order one and exponentially decaying correlations, and no convergence of covariance kernels is needed. The upper bound uses this through localization, and the lower bound uses it through the second Wiener chaos.

\begin{step}[Reduction, Section \ref{sec:red}] \label{step:1}
    By the universality results of Nguyen and Vu \cite{NV21} (Proposition \ref{pro:NV}), the variance of the zero count on the core region $\mathcal S_n$, which consists of four sheets $S_{nj}$ near $\pm1$, differs from its Gaussian counterpart by $o(1)$, and the zeros outside $\mathcal S_n$ contribute $o(\log n)$ to the second moment. It therefore suffices to prove the Gaussian case (Theorem \ref{thm:gauss}). This is the only place in the proof of Theorem \ref{thm:main} where Condition \ref{cond:xi} is used beyond the Gaussian case.
\end{step}

\begin{step}[Geometry, Section \ref{sec:geom}]
    For Gaussian coefficients $\xi=(\xi_0,\dots,\xi_n)$, we write
\[
\frac{P_n(x)}{\sqrt{\Var P_n(x)}}=\langle\xi,g(x)\rangle,
\]
    where $g(x)\in\R^{n+1}$ is the unit vector proportional to $(v_jx^j)_{j=0}^n$. The substitution $x\mapsto1/x$ turns the two outer sheets into inner sheets for the reversed weights $(v_{n-j})_{j=0}^n$, which need not satisfy \eqref{eq:weights} but belong to the class of admissible weights. For admissible weights, the two-sided estimates for the variance function and its derivatives (Lemma \ref{lem:K}) imply that, in the variable $t=\artanh x$, the curve $g$ has speed of order one (Lemma \ref{lem:speed}) and that the correlations of $g$ and of its derivative decay exponentially in $|t-s|$ (Lemma \ref{lem:decay}).
\end{step}

\begin{step}[Upper bound, Section \ref{sec:upper}] \label{step:3}
    We partition $J_n=\artanh(S_{n1})$ into windows of hyperbolic length $1/8$. Jensen's formula bounds the second moment of the count in one window (Lemma \ref{lem:jensen}), so windows at bounded distance have $O(1)$ covariance by the Cauchy--Schwarz inequality. For distant windows, the two-point Kac--Rice formula and the decorrelation estimate \eqref{eq:PsiLip} give covariances that decay exponentially in the distance. Summing yields $\Var N_n(S_{n1})\ll\log n$ (Proposition \ref{prop:upper}), and the Cauchy--Schwarz inequality combines the four sheets.
\end{step}

\begin{step}[Lower bound, Sections \ref{sec:chaos} and \ref{sec:test}] \label{step:4}
    The second Wiener chaos of $\xi$ is the space of centered quadratic forms $\langle\xi\xi^\top-I,\Theta\rangle$ with $\Theta$ symmetric. The projection of $N_n(\mathcal S)$ onto it is $\langle\xi\xi^\top-I,T_{\mathcal S}\rangle$ for an explicit symmetric matrix $T_{\mathcal S}$, obtained from a weighted Rice formula (Lemma \ref{lem:wrice}), so that $\Var N_n(\mathcal S)\ge2\|T_{\mathcal S}\|_{\HS}^2$ (Proposition \ref{prop:chaos}). Such projections are a standard device for variance lower bounds of Gaussian level functionals \cite{GW11, ADL16, BCP19, DNN22}. The new point here is that \eqref{eq:weights} gives no fine information on $g$, so $T_{\mathcal S}$ cannot be computed explicitly. We instead bound $\|T_{\mathcal S}\|_{\HS}$ from below by testing $T_{\mathcal S}$ against order $\log n$ orthogonal vectors supported on disjoint logarithmic scales (Corollary \ref{cor:test}). Each vector is a Mellin-type oscillation, built from sums of the form $\sum_j j^{b-1}e^{-\delta j}$ whose size and phase as $\delta\to0$ are given by the Gamma function. Its frequency is chosen large in terms of $\tau$, $C_1$, and $C_2$ only (Section \ref{sec:test}).
\end{step}

    Theorem \ref{thm:gauss} follows from Steps \ref{step:3} and \ref{step:4}, and Theorem \ref{thm:main} then follows from Step \ref{step:1}. The CLT in Corollary \ref{cor:clt} is a consequence of Theorem \ref{thm:NV}.

    In the non-centered case (Section \ref{sec:noncentered}), Theorem \ref{thm:DNOR} reduces the variance bounds and the CLT for the zeros near $\pm1$ to those of the centered polynomial, while the zeros away from $\pm1$ contribute $O(1)$ to the second moment.

\subsection{Notation and conventions}
    All asymptotic statements are taken in the limit $n\to\infty$, and all implicit constants are independent of $n$. For real quantities $X_n$ and $Y_n$, we write $X_n=O(Y_n)$, $X_n\ll Y_n$, or $Y_n\gg X_n$ if $|X_n|\le C Y_n$ for some constant $C>0$. If both $X_n\ll Y_n$ and $X_n\gg Y_n$ hold, we write $X_n\asymp Y_n$. Finally, we write $X_n\sim Y_n$ if $X_n/Y_n\to1$, and $X_n=o(Y_n)$ if $X_n/Y_n\to0$. Throughout, $C,c>0$ denote constants that may vary from line to line and depend on $\tau,\varepsilon_0,C_0,C_1,C_2$ and any indicated parameters, but not on $n$. The phrase ``for all sufficiently large $n$'' means $n\ge n_0$, with $n_0$ depending only on the same quantities. The coefficients $\xi_j$ are assumed to be i.i.d.\ standard Gaussian in Sections \ref{sec:geom}--\ref{sec:test}; elsewhere, they satisfy Condition \ref{cond:xi}.

\section{Reduction to Gaussian coefficients} \label{sec:red}
    This section reduces Theorem \ref{thm:main} to a statement about Gaussian coefficients on a core region near $\pm1$, and then proves Corollary \ref{cor:clt}. 

    Let $d_n=\exp((\log n)^{1/5})$, $a_n=1/d_n$, and $b_n=d_n/n$. Note that 
\begin{equation}\label{eq:dn}
d_n\ge(\log n)^3\quad\text{for all sufficiently large }n.
\end{equation}
    In particular, $n^Ae^{-d_n/16}\to0$ for every fixed $A$.

    Set
\[
S_{n1}=[1-a_n,1-b_n],\quad S_{n2}=-S_{n1},\quad S_{n3}=\{1/x: x\in S_{n1}\},\quad S_{n4}=-S_{n3},
\]
    and 
\[
\mathcal S_n=S_{n1}\cup S_{n2}\cup S_{n3}\cup S_{n4}.
\]
    We call $S_{n1},\dots,S_{n4}$ the \textit{sheets}, where $S_{n1}$ is a left neighborhood of $1$, $S_{n2}$ is its reflection in the origin, and $S_{n3}$ and $S_{n4}$ are the images of $S_{n1}$ and $S_{n2}$ under $x\mapsto1/x$.

    Let $P_n$ be as in \eqref{eq:poly} under Conditions \ref{cond:xi} and \ref{cond:main}. Let
\[
\widetilde P_n(x):=\sum_{j=0}^n\widetilde\xi_jv_jx^j,
\]
    where $\widetilde\xi_0,\dots,\widetilde\xi_n$ are i.i.d.\ standard Gaussian, and let $\widetilde N_n(S)$ be the number of real zeros of $\widetilde P_n$ in $S$.

    The following proposition is a special case of \cite[Corollary 2.2 and Proposition 2.3]{NV21}.

\begin{proposition}[\cite{NV21}]
\label{pro:NV}
    Let $P_n$ be defined by \eqref{eq:poly}, where Conditions \ref{cond:xi} and \ref{cond:main} hold.  Fix $S\in \{\mathcal S_n, S_{n1}\cup S_{n3}, S_{n2}\cup S_{n4}\}$. There exist constants $C,c>0$, depending only on $\tau$, $\varepsilon_0$, $C_0$, $ C_1$ and $C_2$,  such that for all sufficiently large $n$,
\[
|\Var N_n(S)-\Var\widetilde N_n(S)| \le C (a_n^c+n^{-c}) 
\]
    and
\[
\E N_n^2(\mathbb R\setminus\mathcal S_n) \le C|\log a_n|^4 = C(\log n)^{4/5}.
\]
\end{proposition}

    Let $\R_+=(0,\infty)$ and $\R_-=(-\infty, 0)$. For $S\in \{\mathcal S_n, S_{n1}\cup S_{n3}, S_{n2}\cup S_{n4}\}$, let 
\[
\R_S=\begin{cases}
    \R&\text{if}\quad S=\mathcal S_n,\\
    \R_+&\text{if}\quad S=S_{n1}\cup S_{n3},\\
    \R_-&\text{if}\quad S=S_{n2}\cup S_{n4}.
\end{cases}
\]

\begin{lemma}\label{lem:red}
    Let $S\in \{\mathcal S_n, S_{n1}\cup S_{n3}, S_{n2}\cup S_{n4}\}$. There exists a constant $C>0$ such that for all sufficiently large $n$,
\[
\Bigl|\sqrt{\Var N_n(\R_S)}-\sqrt{\Var\wt N_n(S)}\Bigr|\le C(\log n)^{2/5},
\]
    provided that $\Var\wt N_n(S)\ge 1$. Thus, if $\Var\wt N_n(S) \asymp \log n$, then $\Var N_n(\R_S) \asymp \log n$.
\end{lemma}

\begin{proof}
    Since $\R_S\setminus S\subset\R\setminus\mathcal S_n$, we have 
\[
\Var N_n(\R_S\setminus S)\le\E N_n^2(\R_S\setminus S)\le\E N_n^2(\R\setminus\mathcal S_n).
\]
    The triangle inequality for the standard deviation and Proposition \ref{pro:NV} thus give
\[ 
\left|\sqrt{\Var N_n(\mathbb R_S)}-\sqrt{\Var N_n(S)}\right| \le \sqrt{\Var N_n(\mathbb R_S\setminus S)} \le C(\log n)^{2/5}. 
\] 
    Since $\Var\wt N_n(S)\ge1$ and $C(a_n^c+n^{-c})\le1$ for large $n$, Proposition \ref{pro:NV} yields
\[ 
\Bigl|\sqrt{\Var N_n(S)}-\sqrt{\Var\widetilde N_n(S)}\Bigr|=\frac{\Bigl|\Var N_n(S)-\Var\widetilde N_n(S)\Bigr|}{\sqrt{\Var N_n(S)}+\sqrt{\Var\widetilde N_n(S)}}\le 1. 
\] 
    The first claim follows from the triangle inequality. The second follows because $(\log n)^{2/5}=o(\sqrt{\log n})$.
\end{proof}

    The proof of Theorem \ref{thm:main} thus reduces to the following Gaussian statement, proved in Sections \ref{sec:geom}--\ref{sec:test}.

\begin{theorem}[Gaussian core]\label{thm:gauss}
    For $S\in \{\mathcal S_n, S_{n1}\cup S_{n3}, S_{n2}\cup S_{n4}\}$, there are constants $C,c>0$ such that for all sufficiently large $n$,
\[
c\log n\le\Var\wt N_n(S)\le C\log n.
\]
\end{theorem}

\begin{proof}[Proof of Theorem \ref{thm:main}]
    Apply Lemma \ref{lem:red} and Theorem \ref{thm:gauss} with $S=\mathcal S_n$, for which $\R_S=\R$.
\end{proof}

\begin{proof}[Proof of Corollary \ref{cor:clt}]
    By Theorem \ref{thm:main}, \eqref{eq:lowbound-ass} holds under Conditions \ref{cond:xi} and \ref{cond:main}, and Theorem \ref{thm:NV} gives the CLT.
\end{proof}

    Applying Lemma \ref{lem:red} and Theorem \ref{thm:gauss} with $S\in \{S_{n1}\cup S_{n3}, S_{n2}\cup S_{n4}\}$ gives the following one-sided bound, used in the proof of Theorem \ref{thm:noncentered}.

\begin{proposition}\label{prop:onesided}
    For $\R_S\in\{\R_-, \R_+\}$, there are constants $C,c>0$ such that for all sufficiently large $n$,
\[
c\log n\le\Var N_n(\R_S)\le C\log n.
\]
\end{proposition}

\section{Geometry of the normalized Gaussian process}\label{sec:geom}
    In Sections \ref{sec:geom}--\ref{sec:test}, $\xi=(\xi_0,\ldots,\xi_n)$ is a standard Gaussian vector in $\mathbb R^{n+1}$, and $N_n(S)$ denotes the number of zeros in $S$ of the Gaussian polynomial (written $\widetilde N_n(S)$ in Section \ref{sec:red} and Theorem \ref{thm:gauss}). We work on the core interval
\[
S_{n1}=[1-a_n,1-b_n],\qquad a_n=d_n^{-1},\quad b_n=d_n/n,
\]
    with $d_n=\exp((\log n)^{1/5})$ as in Section \ref{sec:red}, so that \eqref{eq:dn} holds. 

    In this section, we establish uniform estimates for the variance functions of Gaussian polynomials with admissible weights and record the Kac--Rice formulas used in Section \ref{sec:upper}.

\subsection{Admissible weights}
    Let $v=(v_j)_{j=0}^n$ and $v^*=(v_{n-j})_{j=0}^n$. Since $(\xi_j)_{j=0}^n$, $((-1)^j\xi_j)_{j=0}^n$, and $(\xi_{n-j})_{j=0}^n$ have the same law, the number of zeros of $P_n$ on each of the four sheets $S_{n1},\ldots,S_{n4}$ has the law of the number of zeros on $S_{n1}$ of
\[
P_w(x)=\sum_{j=0}^n\xi_jw_jx^j,\qquad w \in \{v, v^*\}.
\]
    Only the marginal law of each sheet count is used here. The dependence among the four sheets is handled by the Cauchy--Schwarz inequality in the proof of the upper bound (Section \ref{sec:upper}) and by the linearity of the second-chaos projection in the proof of the lower bound (Section \ref{sec:chaos}).

    The reversed weights $v^*$ need not satisfy \eqref{eq:weights}, so we work with the following class.

\begin{definition}[Admissible weights]
    Fix $\kappa>-1/2$ and let $s_n\ge (n+1)^{-1}$. A vector $w=(w_0,...,w_n)$ of real weights is \textit{admissible} with parameters $(\kappa, s_n)$ if there exist constants $\kappa_0, c, C>0$ such that
\[
c s_n(1+j)^{2\kappa}\le w_j^2\le C s_n(1+j)^{2\kappa},\qquad 0\le j\le n/2,
\]
    and 
\[
w_j^2\le C(n+1)^{\kappa_0},\qquad n/2<j \le n.
\]
    All constants below that depend on an admissible $w$ depend only on $\kappa,\kappa_0,c,C$, and not on $n$ or $s_n$.
\end{definition}

    By Condition \ref{cond:main}, the weights $v$ and $v^*$ are admissible with parameters $(\tau,1)$ and $(0,(n+1)^{2\tau})$, respectively, and $\kappa_0$ can be taken to be $\max\{2|\tau|,1\}$. For $v^*$, the requirement $s_n\ge(n+1)^{-1}$ holds because $2\tau>-1$.

    For admissible $w$, write $\theta=2\kappa+1>0$ and 
\[ 
k(u):=\sum_{j=0}^n w_j^2u^j. 
\] 
    Thus, $k(x^2)= \Var P_w(x)$. 

    To estimate $k$ and its derivatives, we first record a standard power-series estimate.

\begin{lemma} \label{lem:sum}
For every $a>-1$, there exists $C_a>0$ such that
\[ 
\sum_{j=0}^\infty (1+j)^a u^j \le C_a(1-u)^{-a-1}, \qquad 1/2\le u<1. 
\]
\end{lemma}

\begin{proof} 
    For $a\ge0$ and $t\in[j-1,j]$, $j\ge 1$, 
\[
(1+j)^au^j\le 2^a(1+t)^au^{t}\le 2^a(1+t)^ae^{-t(1-u)}
\]
    where the last inequality follows from $\log(1/u)\ge 1-u$. Integrating over $[j-1,j]$, summing over $j\ge 1$, and adding the $j=0$ term, which is bounded by $(1-u)^{-a-1}$, proves the claim.

    For $-1<a<0$, the sequence $(1+j)^a u^j$ is decreasing in $j$. Since $u^j\le e^{-j(1-u)}$, an integral comparison gives
\[
\sum_{j\ge 0}(1+j)^a u^j \le 1+\int_0^\infty (1+t)^a e^{-t(1-u)}\,dt \le C_a(1-u)^{-a-1},
\]
    which completes the proof.
\end{proof}

\begin{lemma}\label{lem:K}
    Let $w$ be admissible and $E_n=[b_n/8,1/64]$. There exist constants $c,C>0$, depending only on the admissibility parameters, such that, for all sufficiently large $n$, uniformly over $\varepsilon\in E_n$,
\[
cs_n\varepsilon^{-\theta-l}
\le k^{(l)}(1-\varepsilon)
\le Cs_n\varepsilon^{-\theta-l},
\qquad l=0,1,2.
\]
\end{lemma}
\begin{proof}
    For $l\ge0$, let $(j)_l=j(j-1)\cdots(j-l+1)$ denote the falling factorial, and write
\[
k^{(l)}(1-\varepsilon)=\sum_{j=l}^n w_j^2(j)_l(1-\varepsilon)^{j-l}.
\]

    For the lower bound, restrict to $j\in ((4\varepsilon)^{-1}, \varepsilon^{-1}]\cap\mathbb Z$. Since $\varepsilon\ge b_n/8$ and $d_n\to\infty$, this set lies in $[16,n/2]$ for large $n$, and it has at least $c\varepsilon^{-1}$ elements. On this index set, $w_j^2\ge cs_n\varepsilon^{-2\kappa}$, $(j)_l\ge c\varepsilon^{-l}$, and $(1-\varepsilon)^{j-l}\ge(1-\varepsilon)^{1/\varepsilon}\ge c$. Hence, 
\[
k^{(l)}(1-\varepsilon)\ge cs_n\varepsilon^{-2\kappa-l-1}.
\]

    For the upper bound, admissibility and Lemma \ref{lem:sum} yield
\[ 
\sum_{j\le n/2}w_j^2(j)_l(1-\varepsilon)^{j-l} \le Cs_n\sum_{j=0}^\infty(1+j)^{2\kappa+l}(1-\varepsilon)^j \le Cs_n\varepsilon^{-\theta-l}. 
\] 
    For $j>n/2$, $w_j^2\le C(n+1)^{\kappa_0}\le Cs_n(n+1)^{\kappa_0+1}$ because $s_n\ge(n+1)^{-1}$. Using $(1-\varepsilon)^j\le e^{-j\varepsilon}$, $n\varepsilon/2\ge d_n/16$ and $\varepsilon^{-\theta-l}\ge1$,
\[
\sum_{n/2<j\le n}w_j^2(j)_l(1-\varepsilon)^{j-l}\le Cs_n(n+1)^{\kappa_0+l+2}e^{-d_n/16}\le Cs_n\varepsilon^{-\theta-l},
\]
    where the last inequality follows from \eqref{eq:dn}.
\end{proof}

\subsection{The normalized Gaussian process}

    Let $w$ be admissible. For $x\in\R\setminus\{0\}$, define
\[
w(x)=(w_jx^j)_{j=0}^n,\qquad g(x)=\frac{w(x)}{\|w(x)\|},\qquad G(x)=\langle\xi,g(x)\rangle
=\frac{P_w(x)}{\sqrt{k(x^2)}}.
\]
    Then $G$ is a smooth, unit-variance Gaussian process with the same zeros as $P_w$, and its covariance kernel is
\[
r(x,y)=\langle g(x),g(y)\rangle
=\frac{k(xy)}{\sqrt{k(x^2)k(y^2)}}.
\]
    We write $r_{ij}(x,y)=\partial_x^i\partial_y^j r(x,y)$. Since $\|g\|=1$, $g(x)\perp g'(x)$. Moreover, $g'(x)\ne0$ for $x\ne0$ (the vectors $w(x)$ and $w'(x)$ are linearly independent, because the $2\times2$ minor of their zeroth and first coordinates equals $w_0w_1\ne 0$). We may therefore define
\[
h(x)=\frac{g'(x)}{\|g'(x)\|} \qquad\text{and}\qquad  H(x)=\langle\xi,h(x)\rangle.
\]

    For $x\in S_{n1}$, put $\eta(x)=\artanh x$ and $J_n=\eta(S_{n1})=[t_n,u_n]$. Under $x=\tanh t$,
\[
1-x^2=\sech^2t,\qquad dx=(1-x^2)\,dt,\qquad e^{-2t_n}\asymp a_n,\qquad e^{-2u_n}\asymp b_n,
\]
    and hence,
\begin{equation}\label{eq:Jn}
|J_n|\le \frac12\log n\,(1+o(1)), \qquad t_n\to\infty.
\end{equation}

\begin{lemma}[Speed]\label{lem:speed}
    Uniformly for $x\in S_{n1}$,
\[
\|g'(x)\|\asymp(1-x^2)^{-1}.
\]
    Consequently, if $x=\tanh t$ and $\sigma(t):=\|\frac{d}{dt} g(\tanh t)\|$, there are constants $\sigma_1,\sigma_2>0$ such that
\[
\sigma_1 \le \sigma(t) \le \sigma_2, \qquad t\in J_n.
\]
\end{lemma}

\begin{proof}
    Set $\varepsilon=1-x^2$. For $x\in S_{n1}$ and sufficiently large $n$, $\varepsilon\in E_n$. Differentiating $g(x)=w(x)/\|w(x)\|$ gives
\[
\|g'(x)\|^2 = \frac{\|w'(x)\|^2\|w(x)\|^2 - \langle w(x),w'(x)\rangle^2}{\|w(x)\|^4}.
\]
    By Lemma \ref{lem:K},
\[
\|w'(x)\|^2 = k'(x^2)+x^2k''(x^2)\le Cs_n\varepsilon^{-\theta-2} \qquad\text{and}\qquad \|w(x)\|^2 = k(x^2) \ge cs_n\varepsilon^{-\theta}.
\]
    This gives 
\[
\|g'(x)\|^2\le \frac{\|w'(x)\|^2}{\|w(x)\|^2} \le C\varepsilon^{-2}.
\]

    For the lower bound, Lagrange's identity yields
\[
\|w'(x)\|^2\|w(x)\|^2-\langle w(x),w'(x)\rangle^2 = \frac12\sum_{i,j=0}^n
w_i^2w_j^2x^{2(i+j)-2}(i-j)^2.
\]
    Restrict to $(i,j)\in J\times J'$, where
\[
J=\left(\frac1{8\varepsilon},\frac1{4\varepsilon}\right] \cap \mathbb Z,
\qquad J'=\left(\frac1{2\varepsilon},\frac3{4\varepsilon}\right] \cap \mathbb Z.
\]
    Both sets lie in $[0,n/2]$ and have at least $c\varepsilon^{-1}$ elements. For $(i,j)\in J\times J'$,
\[
|i-j|\ge c\varepsilon^{-1}, \qquad x^{2(i+j)-2}\ge c, \qquad w_i^2w_j^2\ge cs_n^2\varepsilon^{-4\kappa}.
\]
    It follows that
\[
\|w'(x)\|^2\|w(x)\|^2-\langle w(x),w'(x)\rangle^2 \ge cs_n^2\varepsilon^{-2\theta-2}.
\]
    Since $\|w(x)\|^4\le Cs_n^2\varepsilon^{-2\theta}$, we obtain $\|g'(x)\|^2\ge c\varepsilon^{-2}$.

    Finally, since $dx/dt=1-x^2$,
\[
\sigma(t)=(1-x^2)\|g'(x)\|\asymp1,
\]
    which proves the second assertion.
\end{proof}

    Since $g(x)\perp h(x)$, the variables $G(x)$ and $H(x)$ are independent standard Gaussians, and $G'(x)=\|g'(x)\|H(x)$. For all $x,y\in\R\setminus\{0\}$, let 
\[
E(x,y):= \begin{pmatrix}
r(x,y)&r_{01}(x,y)/\|g'(y)\|\\
r_{10}(x,y)/\|g'(x)\|&
r_{11}(x,y)/(\|g'(x)\|\|g'(y)\|)
\end{pmatrix}
\]
    be the cross-covariance matrix of $(G(x),H(x))$ and $(G(y),H(y))$. The following lemma bounds it on $S_{n1}\times S_{n1}$.

\begin{lemma}[Decay of correlations]\label{lem:decay}
    Uniformly for $x,y\in S_{n1}$ and $i,j\in\{0,1\}$,
\[
|r_{ij}(x,y)|
\le C\frac{\sech^\theta\bigl(\eta(x)-\eta(y)\bigr)}{(1-x^2)^{i}(1-y^2)^{j}}.
\]
    Consequently, the cross-covariance matrix $E$ satisfies, with $\|\cdot\|$ the operator norm,
\begin{equation}\label{eq:Edecay}
\|E(x,y)\|
\le C e^{-\theta|\eta(x)-\eta(y)|}.
\end{equation}
\end{lemma}

\begin{proof}
    Put $\varepsilon_x=1-x^2$, $\varepsilon_y=1-y^2$, and $\Delta=\eta(x)-\eta(y)$. A direct computation gives
\[
\frac{\varepsilon_x\varepsilon_y}{(1-xy)^2}=\sech^2\Delta.
\]
    Since $b_n\le1-xy\le2a_n$ and $\varepsilon_x,\varepsilon_y\le2(1-xy)$, Lemma \ref{lem:K} applies at $xy,x^2,y^2$ and gives
\[
|r(x,y)| \le C\left(\frac{\varepsilon_x\varepsilon_y}{(1-xy)^2}\right)^{\theta/2} = C\sech^\theta\Delta.
\]
    Next, 
\[
r_{10}(x,y) = \frac{yk'(xy)}{\sqrt{k(x^2)k(y^2)}} -r(x,y)\frac{xk'(x^2)}{k(x^2)}.
\]
    Lemma \ref{lem:K} and $(1-xy)^{-1}\le 2 \varepsilon_x^{-1}$ give
\[
|r_{10}(x,y)| \le C \varepsilon_x^{-1} \sech^\theta\Delta.
\]
    The symmetric estimate holds for $r_{01}$.

    Now, 
\[
r_{11}=\frac{k'(xy)+xyk''(xy)}{\sqrt{k(x^2)k(y^2)}}-\frac{y\,k'(xy)}{\sqrt{k(x^2)k(y^2)}}\,\frac{yk'(y^2)}{k(y^2)}-\frac{xk'(x^2)}{k(x^2)}\,r_{01}.
\]
    Lemma \ref{lem:K} and $(1-xy)^{-2}\le 4(\varepsilon_{x}\varepsilon_y)^{-1}$ give
\[
|r_{11}(x,y)| \le C\varepsilon_x^{-1}\varepsilon_y^{-1} \sech^\theta\Delta.
\]

    Finally, \eqref{eq:Edecay} follows from the estimates for $r_{ij}$, Lemma \ref{lem:speed}, and the inequality
$\sech u\le2e^{-|u|}$.
\end{proof}

\subsection{Kac--Rice formulas}

    Recall that $G(x)$ and $H(x)$ are independent standard Gaussians. The covariance matrix of $(G(x),H(x),G(y),H(y))$ is
\[
\Sigma(E)= \begin{pmatrix}
I_2&E\\
E^\top&I_2
\end{pmatrix}, \qquad E=E(x,y).
\]
    For distinct $x,y\ne0$, let $\Omega(x,y)=\Cov\bigl((H(x),H(y))\mid G(x)=G(y)=0\bigr)$ and
\[
\Psi(x,y)= \frac{\mathbb E|U_1U_2|}{2\pi\sqrt{1-r(x,y)^2}},
\qquad (U_1,U_2)\sim\mathcal N(0,\Omega(x,y)).
\]

\begin{lemma}[Kac--Rice]\label{lem:KR}
    Let $S, S'\subset\mathbb R\setminus\{0\}$ be bounded disjoint intervals, and let $N_w(S)$ be the number of zeros of $P_w$ in $S$. Then
\[
\E N_w(S)=\int_S\rho_1(x)\,dx,
\qquad
\rho_1(x)=\frac{1}{\pi}\|g'(x)\|,
\]
    and
\[
\E[N_w(S)N_w(S')] = \iint_{S\times S'} \rho_2(x,y)\,dx\,dy,
\]
    where
\[
\rho_2(x,y)=\|g'(x)\|\|g'(y)\|\Psi(x,y).
\]
    Moreover, $\Psi(x,y)$ depends on $(x,y)$ only through $E(x,y)$, equals $\pi^{-2}$ when $E=0$, and there is an absolute constant $C$ such that
\begin{equation}\label{eq:PsiLip}
|\Psi(x,y)-\pi^{-2}| \le C\|E(x,y)\|\qquad\text{whenever }\|E(x,y)\|\le 1/2.
\end{equation}
\end{lemma}

\begin{proof}
    Since $G$ has smooth paths and unit variance, and $(G(x),G'(x))$ is nondegenerate, the one-point Kac--Rice formula \cite[Theorem 3.2]{AW09} gives
\[
\rho_1(x)=p_{G(x)}(0)\,\mathbb E\bigl[|G'(x)|\mid G(x)=0\bigr]=\frac1\pi\|g'(x)\|.
\]
    For $x\ne y$, the vectors $g(x)$ and $g(y)$ are linearly independent, since the minor of their zeroth and first coordinates equals $w_0w_1(y-x)/(\|w(x)\|\|w(y)\|)\neq0$. Hence, $(G(x),G(y))$ is nondegenerate for all distinct $x,y\neq0$. Let $I\subset\mathbb R\setminus\{0\}$ be a compact interval containing $S\cup S'$. The second factorial-moment case of \cite[Theorem 3.2]{AW09} on $I$ gives the density
\[
\rho_2(x,y)=p_{G(x),G(y)}(0,0)\,\mathbb E\bigl[|G'(x)G'(y)|\mid G(x)=G(y)=0\bigr],
\]
    and, applying it to $S$, $S'$ and $S\cup S'$ and polarizing, we obtain 
\[
\mathbb E[N_w(S)N_w(S')] = \iint_{S\times S'}\rho_2(x,y)\,dx\,dy.
\]
    This is the stated expression.

    It remains to prove \eqref{eq:PsiLip}. Note that $r=E_{11}$ and that $\Omega=(\Omega_{ij})_{i,j=1}^2$ is the Schur complement in $\Sigma(E)$ of the block of $(G(x),G(y))$. Both depend only on $E$. If $E=0$, the two pairs are independent and 
\[
\Psi=(2\pi)^{-1}(\mathbb E|Z|)^2=\pi^{-2},\qquad Z\sim\mathcal N(0,1).
\]

    If $\|E\|\le1/2$, the eigenvalues of $\Sigma(E)$ lie in $[1/2,3/2]$. Since $\Omega^{-1}$ is the $(H(x),H(y))$ block of $\Sigma(E)^{-1}$, interlacing gives $\frac12 I_2\le\Omega\le\frac32 I_2$. The correlation $\delta=\Omega_{12}/\sqrt{\Omega_{11}\Omega_{22}}$ then satisfies
\[
\delta^2=1-\frac{\det\Omega}{\Omega_{11}\Omega_{22}}\le1-\frac19=\frac89,
\]
    and
\[
\E|U_1U_2|=\frac2\pi\sqrt{\Omega_{11}\Omega_{22}}\bigl(\sqrt{1-\delta^2}+\delta\arcsin\delta\bigr).
\]
    Thus, $\Psi$ is a smooth function of $E$ on the compact convex set $\{\|E\|\le1/2\}$, hence Lipschitz with an absolute constant.
\end{proof}

\section{The upper bound}\label{sec:upper}
    We prove the upper bound in Theorem \ref{thm:gauss}. The argument has two ingredients: a uniform second-moment bound for zero counts on short hyperbolic intervals, obtained from Jensen's formula, and exponential decay of the covariance between counts on well-separated intervals, obtained from the two-point Kac--Rice formula.

\begin{lemma}[Local bounds]\label{lem:jensen}
    Let $w$ be admissible, $t_0\in J_n$, and let $S\subset S_{n1}$ be an interval such that $\eta(S)\subset[t_0,t_0+1/8]$. Then there is a constant $C$, independent of $n$, $t_0$, and $S$, such that for all sufficiently large $n$,
\[
\mathbb E N_w(S)^2\le C.
\]
\end{lemma}

\begin{proof}
    Put $x_0=\tanh t_0$ and $\varepsilon=1-x_0$, so $b_n\le\varepsilon\le a_n$. Since
\[
\tanh(t_0+1/8)-\tanh t_0 \le \sinh(1/8)(1-x_0^2)\le 2\sinh(1/8)\,\varepsilon < 0.26\varepsilon,
\]
    the interval $S$ is contained in the disk $D(x_0,\varepsilon/2)$. Let 
\[
\mathcal M=\max_{|z-x_0|=3\varepsilon/4}|P_w(z)|.
\]
    Since $P_w(x_0)\ne0$ almost surely, Jensen's formula \cite[Theorem 15.18]{R87} on the disk $D(x_0,3\varepsilon/4)$ gives
\[
N_w(S)\log(3/2)\le\log\mathcal M-\log|P_w(x_0)|.
\]
    By subharmonicity of $|P_w|^2$, for $|z-x_0|=3\varepsilon/4$,
\[
|P_w(z)|^2\le\frac{64}{\pi\varepsilon^2}\int_{D(x_0,7\varepsilon/8)}|P_w(\zeta)|^2\,dA(\zeta).
\]
    Since $D(x_0,7\varepsilon/8)\subset\{|\zeta|\le1-\varepsilon/8\}$ and $k$ is increasing, $\E|P_w(\zeta)|^2\le k(1-\varepsilon/8)$ on this disk. Hence, 
\[
\E\mathcal M^2\le49\,k(1-\varepsilon/8)\le Cs_n\varepsilon^{-\theta}
\]
    by Lemma \ref{lem:K}. Since $k(x_0^2)\ge cs_n\varepsilon^{-\theta}$, the variable $Y=\mathcal M/\sqrt{k(x_0^2)}$ satisfies $\E Y^2\le C$. Writing $P_w(x_0)=\sqrt{k(x_0^2)}\,Z$ with $Z\sim\mathcal N(0,1)$, we obtain
\[
N_w(S)\le\frac{\log_+Y+\log_+(1/|Z|)}{\log(3/2)},\qquad \log_+(u)=\max\{0, \log u\}.
\]
    Since $(\log_+Y)^2\le Y^2$ and $\E(\log_+(1/|Z|))^2<\infty$, it follows that $\E N_w(S)^2\le C$.
\end{proof}

\begin{proposition}[Upper bounds]\label{prop:upper}
    Let $w$ be admissible. Then, for all sufficiently large $n$ and with $C$ depending only on the admissibility constants,
\[
\mathbb EN_w(S_{n1})\le C|J_n| = O(\log n), \qquad \Var N_w(S_{n1})\le C|J_n| = O(\log n).
\]
    Consequently, for every union $\mathcal S$ of some of the four sheets $S_{n1},\ldots,S_{n4}$,
\[
\E N_n(\mathcal S)=O(\log n), \qquad \Var N_n(\mathcal S)=O(\log n).
\]
\end{proposition}

\begin{proof}
    We first consider the expectation. By the one-point Kac--Rice formula, Lemma \ref{lem:speed}, and \eqref{eq:Jn},
\[
\mathbb EN_w(S_{n1}) = \frac1\pi\int_{S_{n1}}\|g'(x)\|\,dx = \frac1\pi\int_{J_n}\sigma(t)\,dt \le \frac{\sigma_2}{\pi}|J_n| = O(\log n).
\]

    For the variance, choose $j_0\ge1$, depending only on the fixed parameters, such that $Ce^{-\theta j_0}\le 1/2$, where $C$ is the constant in \eqref{eq:Edecay}. Partition $J_n$ into the intervals
\[
T_j=[t_n+(j-1)/8,\;t_n+j/8)\cap J_n,\qquad 1\le j\le j_n:=\lceil8|J_n|\rceil,
\]
    and let $S_j=\tanh(T_j)$. Since a Gaussian polynomial almost surely does not vanish at a prescribed point, $N_w(S_{n1})=\sum_{j=1}^{j_n}N_w(S_j)$ almost surely, and
\[
\Var N_w(S_{n1}) = \sum_{j,l=1}^{j_n} \Cov(N_w(S_j),N_w(S_l)).
\]
    We separate the pairs according to their hyperbolic distance.

    \textit{Near pairs.} If $|j-l|\le 8j_0$, the Cauchy--Schwarz inequality and Lemma \ref{lem:jensen} give 
\[|\Cov(N_w(S_j),N_w(S_l))|\ll 1.\]
    For each $j$, there are at most $16j_0+1$ such $l$.

    \textit{Far pairs.} If $|j-l|\ge8j_0+1$, the hyperbolic distance between $S_j$ and $S_l$ is at least $(|j-l|-1)/8\ge j_0$. The intervals are then disjoint, and Lemma \ref{lem:KR} gives
\[
\Cov(N_w(S_j),N_w(S_l)) =\iint_{S_j\times S_l} \bigl(\rho_2(x,y)-\rho_1(x)\rho_1(y)\bigr)\,dx\,dy.
\]
    By \eqref{eq:Edecay} and the choice of $j_0$, $\|E(x,y)\|\le1/2$ on $S_j\times S_l$, so \eqref{eq:PsiLip} yields
\begin{align*}
|\rho_2(x,y)-\rho_1(x)\rho_1(y)| &\ll \|g'(x)\|\|g'(y)\|\|E(x,y)\|\\
& \le C\|g'(x)\|\|g'(y)\| e^{-\theta|\eta(x)-\eta(y)|}.
\end{align*}
    Using Lemma \ref{lem:speed} and $dx=(1-x^2)\,dt$,
\[
|\Cov(N_w(S_j),N_w(S_l))| \le C\int_{T_j}\int_{T_l} e^{-\theta|t-s|}\, dtds \ll e^{-\theta(|j-l|-1)/8},
\]
    and, uniformly in $j$,
\[
\sum_{l:\,|j-l|\ge8j_0+1}e^{-\theta(|j-l|-1)/8}\ll 1.
\]
    Summing the near and far contributions, we obtain 
\[
\Var N_w(S_{n1}) \ll j_n\le C|J_n|.
\]

    Finally, by the reduction in Section \ref{sec:geom}, each sheet count has the law of $N_w(S_{n1})$ for some admissible $w$, so its mean and variance are $O(\log n)$. For a union $\mathcal S$ of $m\le4$ sheets, summation gives the bound on $\E N_n(\mathcal S)$, and the Cauchy--Schwarz inequality gives 
\[
\Var N_n(\mathcal S)\le m\sum_{S_{nj}\subset\mathcal S}\Var N_n(S_{nj})=O(\log n).
\]
\end{proof}

\section{The second Wiener chaos}\label{sec:chaos} 
    The lower bound in Theorem \ref{thm:gauss} comes from the second Wiener chaos of the Gaussian coefficient vector. This section identifies the projection of the number of real zeros onto this chaos and reduces the lower bound to an estimate for a quadratic form.

    Throughout, $g$ and $h$ are as in Section \ref{sec:geom} with the original weights $w=v$ on $\R\setminus\{0\}$, i.e., 
\[
g(x)=\frac{v(x)}{\|v(x)\|}\qquad\text{and}\qquad h(x)=\frac{g'(x)}{\|g'(x)\|}.
\]
    To describe the geometry on the outer sheets, let $z=1/x$ and
\[
k^*(u)=\sum_{j=0}^n v_{n-j}^2u^j, \qquad g^*(z)= \biggl(\frac{v_{n-j}z^j}{\sqrt{k^*(z^2)}}\biggr)_{j=0}^n.
\]
    Since $k(x^2)=x^{2n}k^*(z^2)$,
\begin{equation}\label{eq:recip}
g(x)=\sgn(x)^n\mathcal R g^*(z), \qquad z=1/x,
\end{equation}
    where $\mathcal R$ is the coordinate reversal. Thus, up to an isometry and a change of parameter, the curve $g$ on $S_{n3}\cup S_{n4}$ is the curve of the admissible weights $v^*$ on the inner sheets. In particular, $g'\ne0$ on all four sheets, and by Lemma \ref{lem:speed}, the arclength $ds=\|g'(x)\|\,|dx|=\|g^{*\prime}(z)\|\,|dz|$ satisfies $ds\le C\,dt$ with $t=\artanh|z|$.

    For a union $\mathcal S$ of some of the sheets, let $l(\mathcal S)$ be the set of $j$ with $S_{nj}\subset\mathcal S$, and define the symmetric matrix
\begin{equation}
    \label{eq:TS}
    T_{\mathcal S}=\frac1{2\pi}\sum_{j\in l(\mathcal S)}\int_{S_{nj}}\bigl( h h^\top-gg^\top\bigr)\,ds,\qquad ds=\|g'(x)\|\,|dx|.
\end{equation}
    The integrand does not depend on the orientation of $h$. 

\begin{lemma}[Weighted Rice formula]\label{lem:wrice}
    For every symmetric $(n+1)\times(n+1)$ matrix $\Theta$ and every $j$,
\[
\E\bigl[N_n(S_{nj})\,\langle\xi\xi^\top-I,\Theta\rangle\bigr] = \frac1\pi\int_{S_{nj}}\bigl\langle\Theta, h h^\top-gg^\top\bigr\rangle\,ds.
\]
\end{lemma}

\begin{proof}
    Let $F(\xi)=\langle\xi\xi^\top-I,\Theta\rangle$ and, for $D>0$, let $F_D=\max\{-D,\min\{D,F\}\}$. Fix $n+1$ distinct points $x_0,\dots,x_n\in S_{nj}$. The vectors $g(x_i)$ are linearly independent, so $\xi=A\,(G(x_i))_{i=0}^n$ for a deterministic matrix $A$. Hence, $F_D(\xi)$ is a bounded continuous function of the sample path of $G$ on $S_{nj}$. Since $G$ has smooth paths and $(G(x),G'(x))$ is nondegenerate for each $x\in S_{nj}$, the weighted Rice formula \cite[Theorem 6.4]{AW09} applies and gives
\[
\mathbb E \bigl[N_n(S_{nj})F_D(\xi)\bigr] = \frac1{\sqrt{2\pi}}\int_{S_{nj}} \mathbb E \bigl[|G'(x)|\,F_D(\xi)\mid G(x)=0\bigr]\,dx.
\]
    Let $D\to\infty$. As $N_n(S_{nj})\le n$ almost surely and $\E|F(\xi)|<\infty$, the left-hand side converges to $\E[N_n(S_{nj})F(\xi)]$ by dominated convergence. 

    For the right-hand side, conditionally on $G(x)=0$, $\xi$ is standard Gaussian on $g(x)^\perp$. Writing $P=I-gg^\top$, we have
\[
\E[F(\xi)^2\mid G(x)=0] = 2\|P\Theta P\|_{\HS}^2+(g^\top\Theta g)^2 \le 3\|\Theta\|_{\HS}^2.
\]
    Since $G'(x)=\|g'(x)\|\langle\xi,h(x)\rangle$, Cauchy--Schwarz gives
\[
\E[|G'(x)|\,|F_D(\xi)|\mid G(x)=0] \le \sqrt{3}\,\|\Theta\|_{\HS}\|g'(x)\|,
\]
    which is integrable over the compact interval $S_{nj}\subset\R\setminus\{0\}$. Dominated convergence gives
\[
\E[N_n(S_{nj})F(\xi)] = \frac1{\sqrt{2\pi}}\int_{S_{nj}} \E[|G'(x)|F(\xi)\mid G(x)=0]\,dx.
\]

    To evaluate the integrand, decompose
\[
\xi=Hh+\zeta,
\]
    where, conditionally on $G(x)=0$, $H\sim\mathcal N(0,1)$ is independent of $\zeta$, and
\[
\E[\zeta\zeta^\top] = I-gg^\top-hh^\top.
\]
    Using $|G'(x)|=\|g'(x)\||H|$, $\E|H|=\sqrt{2/\pi}$, and $\E|H|^3=2\sqrt{2/\pi}$, we obtain
\begin{align*}
\E[|H|(\xi\xi^\top-I)\mid G(x)=0]
    &= \E|H|^3\,hh^\top + \E|H|(I-gg^\top-hh^\top) -\E|H|I\\
    &= \sqrt{\frac2\pi}(hh^\top-gg^\top).
\end{align*}
    Substituting into the Rice integral proves the lemma.
\end{proof}

\begin{proposition}\label{prop:chaos}
    Let $\mathcal S$ be a union of some of the sheets and let $\Pi_2 N_n(\mathcal S)$ denote the projection of $N_n(\mathcal S)$ onto the second Wiener chaos of $\xi$. Then 
\[
\Pi_2N_n(\mathcal S)=\langle\xi\xi^\top-I,T_{\mathcal S}\rangle,
\]
    where $T_{\mathcal S}$ is given by \eqref{eq:TS}. Consequently, 
\[
\Var N_n(\mathcal S) \ge \Var \Pi_2N_n(\mathcal S) =2\|T_{\mathcal S}\|_{\HS}^2.
\]
\end{proposition}

\begin{proof}
    The second Wiener chaos consists of the variables $\langle\xi\xi^\top-I,\Theta\rangle$, with $\Theta$ symmetric, and 
\[
\E\bigl[\langle\xi\xi^\top-I,\Theta\rangle\langle\xi\xi^\top-I,\Theta'\rangle\bigr]=2\langle\Theta,\Theta'\rangle_{\HS}.
\]
    Hence, if $\Pi_2N_n(\mathcal S) =\langle\xi\xi^\top-I,\Theta_{\mathcal S}\rangle$,
    then
\[
2\langle\Theta_{\mathcal S},\Theta\rangle_{\HS} = \E[N_n(\mathcal S) \langle\xi\xi^\top-I,\Theta\rangle]
\]
    for every symmetric $\Theta$. Summing Lemma \ref{lem:wrice} over the sheets in $\mathcal S$ gives
\[
\E[N_n(\mathcal S) \langle\xi\xi^\top-I,\Theta\rangle] = 2\langle T_{\mathcal S},\Theta\rangle_{\HS}.
\]
    Thus, $\Theta_{\mathcal S}=T_{\mathcal S}$, and the variance bound follows from orthogonality of the Wiener chaoses.
\end{proof}

    Since $\Pi_2$ is linear, $T_{\mathcal S}$ is the sum of the matrices of the individual sheets. Hence, no information on the dependence among the sheet counts is needed.

    For $\nu\in\R^{n+1}$, put $\varphi_\nu(x)=\langle\nu,g(x)\rangle$. Since $h=g'/\|g'\|$, 
\[
\frac{d\varphi_\nu}{ds}=\frac{1}{\|g'(x)\|} \frac{d}{dx}\langle\nu,g(x)\rangle=\langle\nu, h(x)\rangle.
\]
    Therefore,
\begin{equation}\label{eq:quad}
\langle\nu,T_{\mathcal S}\nu\rangle=\frac1{2\pi}\sum_{j\in l(\mathcal S)} \int_{S_{nj}}\biggl[\Bigl(\frac{d\varphi_\nu}{ds}\Bigr)^2-\varphi_\nu^2\biggr]\,ds.
\end{equation}

\begin{corollary}\label{cor:test}
    Let $\nu^{(1)},\ldots,\nu^{(m)}\in\mathbb R^{n+1}$ be pairwise orthogonal nonzero vectors such that $\langle\nu^{(j)},T_{\mathcal S}\nu^{(j)}\rangle \ge c\|\nu^{(j)}\|^2$ for some $c>0$ and all $1\le j\le m$. Then
\[
\|T_{\mathcal S}\|_{\HS}^2\ge c^2m \qquad\text{and}\qquad \Var N_n(\mathcal S)\ge 2c^2m.
\]
\end{corollary}

\begin{proof}
    Let $\widehat\nu_j=\nu^{(j)}/\|\nu^{(j)}\|$ and $B=\sum_{j=1}^m\widehat\nu_j\widehat\nu_j^\top$, an orthogonal projection of rank $m$, so $\|B\|_{\HS}=\sqrt m$. Then $\langle T_{\mathcal S},B\rangle_{\HS}=\sum_{j=1}^m\langle\widehat\nu_j,T_{\mathcal S}\widehat\nu_j\rangle\ge cm$, and the Cauchy--Schwarz inequality gives $\|T_{\mathcal S}\|_{\HS}\ge c\sqrt m$. The variance bound follows from Proposition \ref{prop:chaos}.
\end{proof}

\section{Mellin-modulated test vectors}\label{sec:test}
    By Corollary \ref{cor:test}, the lower bound in Theorem \ref{thm:gauss} follows once we construct $m\ge c\log n$ pairwise orthogonal vectors $\nu$ with $\langle\nu,T_{\mathcal S}\nu\rangle\ge c\|\nu\|^2$. By \eqref{eq:quad}, this requires the derivative energy of $\varphi_\nu=\langle\nu,g\rangle$ along the curve to exceed its $L^2$ energy by a margin proportional to $\|\nu\|^2$. Under \eqref{eq:weights}, neither the speed of $g$ nor the normalization $\sqrt{k(x^2)}$ is explicit, so $\varphi_\nu$ must oscillate at a frequency larger than constants that are not known explicitly.

    We take $\nu_0=0$ and $\nu_j=\Real f(j)/v_j$ for $1\le j\le n$, where $f(u)=u^{b-1}W(\log u)$, $\Real b=\theta/2$, $\Img b=\omega/2$, and $W$ is a smooth cutoff of width $2R$ in $\log u$ centered at a point $\lambda_c$ of the logarithmic scale. Dividing by $v_j$ removes the weights from $\varphi_\nu$. On $S_{n1}$, $\varphi_\nu$ is then a Mellin-type sum with amplitude of order $1$ and phase $\frac\omega2\log(1/\delta)$ (Lemma \ref{lem:mellin}). For $\omega=\omega_0$ large and $R\ge R_0$, the sheet $S_{n1}$ contributes at least $c_2R$ to \eqref{eq:quad}, whereas the other sheets cost at most $c_2R/2$ (Lemmas \ref{lem:margin} and \ref{lem:others}). The parameters are fixed in the order $\omega_0$, then $R_0$. Finally, we choose $c\log n$ centers $\lambda_c$ spaced $2R+2$ apart, so that the corresponding vectors have disjoint supports (Section \ref{subsec:windows}).

    Throughout this section, $\xi$ is a standard Gaussian vector, $G(x)=\langle\xi,g(x)\rangle$ is the Gaussian field of Section \ref{sec:geom}, and $N_n(S)$ counts the zeros in $S$ of the Gaussian polynomial $P_n$. Moreover, $w=v$ is admissible with parameters $(\tau,1)$, so $\theta=2\tau+1>0$. Recall $J_n=[t_n,u_n]$ and the constant $\sigma_2$ of Lemma \ref{lem:speed}. Unless indicated otherwise, $\|\cdot\|_2$ is the norm of $L^2(J_n,dt)$.

\subsection{Hyperbolic variables and test vectors} 
    For $t>0$, set
\[
x=\tanh t,\qquad \delta=-\log x=2\artanh(e^{-2t}),\qquad L=-\log\delta.
\]
    Note that functions of $L$ are regarded as functions of $t$ through $L=L(t)$. A direct computation gives $\delta\asymp 1-x^2$, $L_t=2+O(e^{-4t})$, and $L_{tt}=O(e^{-4t})$ as $t\to\infty$. Hence, there is $t_0>0$ such that, for $t\ge t_0$,
\begin{equation}\label{eq:Lt}
\frac32\le L_t\le\frac52,\qquad |L_t-2|+|L_{tt}|\le Ce^{-4t}\le 1.
\end{equation}
    Since $t_n\to\infty$,  $J_n\subset[t_0,\infty)$ for large $n$. At $t_n$, we have $\delta=-\log(1-a_n)\sim a_n$, and at $u_n$, we have $\delta\sim b_n=d_n/n$. This gives
\[
L(t_n)=\log d_n+o(1) \qquad\text{and}\qquad L(u_n)=\log n-\log d_n+o(1).
\]
    Since $\log d_n=(\log n)^{1/5}=o(\log n)$,
\begin{equation}\label{eq:Ldiff}
L(u_n)-L(t_n)\ge\frac12\log n\qquad\text{for all sufficiently large }n.
\end{equation}

    For $x\in S_{n1}$, $b_n\le 1-x^2\le 2a_n\le 1/64$ for large $n$, so Lemma \ref{lem:K} applies with $\varepsilon=1-x^2$, $s_n=1$, $l=0,1$. Since $\frac{d}{dt}k(x^2)=2x(1-x^2)k'(x^2)$, there are constants $0<M_1\le M_2$ and $C_3$, depending only on $\tau,C_1,C_2$, such that for $t\in J_n$,
\begin{equation}\label{eq:m}
M_1\le M(t):=\frac{\delta^{-\theta/2}}{\sqrt{k(x^2)}} \le M_2, \qquad |(\log M)_t|+|(\log k(x^2))_t|\le C_3.
\end{equation}

    Fix $W_0\in C_c^\infty([-1,1])$ satisfying
\[
0\le W_0\le1,\qquad W_0=1\ \text{on }[-1/2, 1/2],\qquad \|W_0'\|_\infty\le 4.
\]
    For $\omega,R\ge1$ and $\lambda_c\in\R$, put
\[
W(\lambda)=W_0\Bigl(\frac{\lambda-\lambda_c}{R}\Bigr),\qquad b=\frac{\theta+i\omega}{2},\qquad f(u)=u^{b-1}W(\log u)\quad(u\ge1).
\]
    Define the real vector $\nu\in\R^{n+1}$ by
\[
\nu_0=0,\qquad \nu_j=\frac{\Real f(j)}{v_j},\quad 1\le j\le n.
\]
    We call $\lambda_c$ an \textit{interior center} if
\[
L(t_n)+1\le\lambda_c-R,\qquad \lambda_c+R\le L(u_n)-1.
\]
    Then $\lambda_c-R>0$ for large $n$, and $W(L(t))\ne0$ only if $L(t)\in[L(t_n)+1,L(u_n)-1]$. Since $L$ is increasing, $W(L(\cdot))$ is supported in $J_n$. Moreover, for large $n$,
\[
e^{\lambda_c+R}\le e^{L(u_n)-1}\le\frac{n}{2d_n}\le\frac n2,
\]
    so $\nu$ is supported in $\{j:e^{\lambda_c-R}\le j\le e^{\lambda_c+R}\}\subset[1,n/2]$. Since $|f(j)|^2\le j^{\theta-2}$ and $v_j^2\ge c\,j^{2\tau}$ with $c=C_1^2\min\{1,2^{2\tau}\}$,
\begin{equation}\label{eq:etanorm}
\|\nu\|^2\le c^{-1}\sum_{e^{\lambda_c-R}\le j\le e^{\lambda_c+R}}\frac1j\le c^{-1}(1+2R)\le C_4R.
\end{equation}

\subsection{Mellin asymptotics}
    We compute $\varphi_\nu$ on $S_{n1}$ and show that the alternating sums governing $S_{n2}$ are negligible. 

    For $a\in\{0,1\}$, let 
\[
S_a(\delta)=\sum_{j\ge1}j^a f(j)e^{-\delta j},
\]
    so that $S_0'=-S_1$. For $x=e^{-\delta}\in(0,1)$, the definition of $\nu$ gives
\begin{equation}\label{eq:phiS}
\varphi_\nu(x)=\frac{1}{\sqrt{k(x^2)}}\Real S_0(\delta),\qquad \varphi_\nu(-x)=\frac{1}{\sqrt{k(x^2)}}\,\Real\sum_{j\ge1}(-1)^jf(j)x^j.
\end{equation}

\begin{lemma}[Mellin asymptotics]\label{lem:mellin}
    Let $\omega,R\ge1$ and let $\lambda_c$ be an interior center. For $t\in J_n$, put $x=\tanh t$ and $\psi_c(t)= S_0(\delta)/\sqrt{k(x^2)}$. There is a constant $C_\omega$, depending only on $\tau,C_1,C_2,\omega$ and not on $R$, $\lambda_c$, or $n$, such that for all sufficiently large $n$,
\begin{align}
\psi_c&=M\,e^{i\omega L/2}\bigl[\Gamma(b)W(L)+r_0\bigr],\label{eq:psi0}\\
\partial_t\psi_c&=M\,e^{i\omega L/2}\bigl[(i\omega+\phi)\Gamma(b)W(L)+r_1\bigr],\label{eq:psi1}
\end{align}
    where $L=L(t)$, $\phi=\theta-\frac12(\log k(x^2))_t$ is real and bounded on $J_n$, and 
\[
\|r_0\|_2^2+\|r_1\|_2^2\le C_\omega.
\]
    Moreover, there is a function $\bar r\ge0$ such that
\begin{equation}\label{eq:alt}
\Bigl|\sum_{j\ge1}(-1)^jf(j)x^j\Bigr|\le M\sqrt{k(x^2)}\,\bar r(L),\qquad \|\bar r\|_2^2\le C_\omega e^{-2(\lambda_c-R)}.
\end{equation}
\end{lemma}

\begin{proof}
    Let $F_a(u)=u^{b+a-1}W(\log u)e^{-\delta u}=u^af(u)e^{-\delta u}$. Then $F_a$ is smooth and supported in $[e^{\lambda_c-R},e^{\lambda_c+R}]\subset(1,\infty)$, and $S_a(\delta)=\sum_{j\in\mathbb Z}F_a(j)$. 

    For $F\in C^1_c(\R)$,
\begin{equation}\label{eq:sumint}
\Bigl|\sum_{j\in \Z}F(j)-\int_\R F(u)\,du \Bigr|\le\int_\R|F'(u)|du,\qquad \Bigl|\sum_{j\in \Z}(-1)^jF(j)\Bigr|\le\int_\R|F'(u)|du.
\end{equation}
    Indeed, the first follows by comparing $F(j)$ with the mean of $F$ on $[j,j+1]$, and the second by pairing $F(2j)-F(2j+1)$.

    \textit{Main term.} Substituting $u=e^{L+y}$, recalling $\delta=e^{-L}$, and using
\[
\int_\R e^{(b+a)y}\exp(-e^y)\,dy=\Gamma(b+a),
\]
    we obtain
\[
\int_\R F_a(u)\,du = e^{(b+a)L}\bigl[\Gamma(b+a)W(L)+I_a(L)\bigr]
\]
    and 
\[
I_a(L) = \int_\R e^{(b+a)y}\bigl[W(L+y)-W(L)\bigr]\exp(-e^y)\,dy.
\]
    Since $|W(L+y)-W(L)|\le\min\{1,4|y|/R\}$, we have $|I_a(L)|\le C/R$ for all $L$. Let $\ell=\dist(L,[\lambda_c-R,\lambda_c+R])$. If $\ell>1$, then $W(L)=0$ and $W(L+y)\neq0$ only for $|y|\ge\ell$, so $|I_a(L)|\le Ce^{-\gamma\ell}$, where $\gamma:=\theta/2+a>0$. Hence,
\[
\|I_a\|_{L^2(\R,dL)}^2\le C\,\frac{R+1}{R^2}+C\ll 1.
\]

    \textit{Error terms.} On the support of $F_a$, using $|W'|\le4/R$ and $|b+a-1|\le|b|+1$,
\[
|F_a'(u)| \le  C(1+\omega)\,u^{\gamma-1}\bigl(u^{-1}+\delta\bigr)e^{-\delta u}.
\]
    Substituting $u=e^{L+y}$ gives
\[
\int_\R|F_a'(u)|du \le C(1+\omega)\,e^{\gamma L}e^{-(\lambda_c-R)}\,\Lambda_a\bigl(L-\lambda_c+R\bigr),
\]
    where
\begin{align*}
\Lambda_a(s)&=e^{-s}\int_{-s}^\infty\bigl(e^{(\gamma-1)y}+e^{\gamma y}\bigr)\exp(-e^y)\,dy\\
&\le \begin{cases}
        C(1+s)e^{-\min\{1,\gamma\}s}&\text{for}\quad s\ge 0,\\
        Ce^{|s|-e^{|s|}/2}&\text{for}\quad s<0.
    \end{cases}
\end{align*}
    Thus, $\|\Lambda_a\|_{L^2(\R)}\le C$, with $C$ depending only on $\theta$.

    \textit{Proofs of \eqref{eq:psi0} and \eqref{eq:alt}.} By \eqref{eq:sumint},
\[
S_a(\delta)=e^{(b+a)L}\bigl[\Gamma(b+a)W(L)+q_a(L)\bigr],
\]
    where 
\[
 q_a=I_a+e^{-(b+a)L}\Bigl(S_a(\delta)-\int_\R F_a(u)du\Bigr),\qquad |e^{-(b+a)L}|=e^{-\gamma L}.
\]
    The estimates above and $\lambda_c-R>0$ give $\|q_a\|_{L^2(dL)}\le C_\omega$. Since $dt=dL/L_t$ and $L_t\ge 3/2$, $L^2(J_n,dt)$-norms are bounded by the corresponding $L^2(dL)$-norms. As $e^{bL}/\sqrt{k(x^2)}=Me^{i\omega L/2}$, \eqref{eq:psi0} holds with $r_0=q_0$. By the second inequality in \eqref{eq:sumint},
\[
\Bigl|\sum_{j\ge1}(-1)^jf(j)x^j\Bigr| \le \int_\R|F_0'(u)|du \le e^{\theta L/2}\bar r(L),
\]
    where 
\[
\bar r(L) := C(1+\omega)e^{-(\lambda_c-R)}\Lambda_0\bigl(L-\lambda_c+R\bigr).
\]
    Since $e^{\theta L/2}=\delta^{-\theta/2}=M\sqrt{k(x^2)}$, this implies \eqref{eq:alt}.

    \textit{Proof of \eqref{eq:psi1}.} From $\delta_t=-\delta L_t$ and $S_0'=-S_1$,
\[
\partial_t\psi_c = \frac{\delta L_tS_1(\delta)}{\sqrt{k(x^2)}}-\frac12\psi_c\,(\log k(x^2))_t.
\]
    Insert $S_1=e^{(b+1)L}[b\Gamma(b)W(L)+q_1]$, using $\Gamma(b+1)=b\Gamma(b)$ and $\delta e^{(b+1)L}=e^{bL}$. Writing $L_tb=2b+(L_t-2)b$ with $2b=\theta+i\omega$, we obtain \eqref{eq:psi1} with
\[
r_1 = (L_t-2)\,b\Gamma(b)W(L)+L_tq_1-\frac12(\log k(x^2))_t\,q_0.
\]
    By \eqref{eq:Lt}, 
\[
\int_{J_n}|L_t-2|^2dt\le C\int_{t_0}^\infty e^{-8t}\,dt\le C.
\]
    The other two terms are controlled by \eqref{eq:m} and the bounds on $q_0,q_1$. Finally, $\phi$ is real and bounded by \eqref{eq:m}.
\end{proof}

\subsection{The quadratic-form margin}
    The oscillation of $\nu$ makes the quadratic form positive on $S_{n1}$, and the other sheets contribute only a small error.

\begin{lemma}[Margin on $S_{n1}$]\label{lem:margin}
    There are constants $\omega_0\ge 1$, $R_0\ge1$, and $c_2>0$, depending only on $\tau,C_1,C_2$, such that for $\omega=\omega_0$, every $R\ge R_0$, every interior center $\lambda_c$, and all sufficiently large $n$,
\[
\frac1{2\pi}\int_{S_{n1}}\biggl[\Bigl(\frac{d\varphi_\nu}{ds}\Bigr)^2-\varphi_\nu^2\biggr]ds \ge c_2R.
\]
\end{lemma}

\begin{proof}
    Let $\psi(t)=\varphi_\nu(\tanh t)=\Real\psi_c(t)$ by \eqref{eq:phiS}. Since $d\varphi_\nu/ds=\psi_t/\sigma$ and $ds=\sigma dt$, Lemma \ref{lem:speed} gives
\begin{equation}\label{eq:basic}
\int_{S_{n1}}\biggl[\Bigl(\frac{d\varphi_\nu}{ds}\Bigr)^2-\varphi_\nu^2\biggr]ds = \int_{J_n}\Bigl(\frac{\psi_t^2}{\sigma}-\sigma\psi^2\Bigr)dt \ge\frac1{\sigma_2}\|\psi_t\|_2^2-\sigma_2\|\psi\|_2^2.
\end{equation}
    Set $A=M\,|\Gamma(b)|\,W(L)$ and $\vartheta=\omega L/2+\arg\Gamma(b)$. Then $A$ is smooth with compact support in $J_n$. Taking real parts in Lemma \ref{lem:mellin},
\[
\psi = A\cos\vartheta+e_0,\qquad \psi_t=-\omega A\sin\vartheta+e_1,
\]
    where $|e_0|\le M_2|r_0|$ and $|e_1|\le\|\phi\|_\infty A+M_2|r_1|$. Hence,
\begin{equation}\label{eq:err}
\|e_0\|_2^2\le M_2^2C_\omega,\qquad \|e_1\|_2^2\le2\|\phi\|_\infty^2\|A\|_2^2+2M_2^2C_\omega.
\end{equation}

    \textit{Size of $A$.} The set where $W(L)=1$ has $t$-length at least $2R/5$ by \eqref{eq:Lt}, so
\begin{equation}\label{eq:Q}
    \|A\|_2^2\ge\frac25\,M_1^2|\Gamma(b)|^2R.
\end{equation}

    \textit{Oscillation.} By \eqref{eq:Lt}, $\vartheta_t=\omega L_t/2\ge3\omega/4$ and $|\vartheta_{tt}|/(2\vartheta_t^2)=|L_{tt}|/(\omega L_t^2)\le1/\omega$. Integrating by parts (there are no boundary terms),
\[
\Bigl|\int_{J_n}A^2\cos2\vartheta\,dt\Bigr| \le \int_{J_n}\frac{|(A^2)_t|}{2\vartheta_t}\,dt + \int_{J_n}\frac{A^2|\vartheta_{tt}|}{2\vartheta_t^2}\,dt \le \frac{2}{3\omega}\int_{J_n}|(A^2)_t|\,dt+\frac{\|A\|_2^2}{\omega}.
\]
    By \eqref{eq:m} and $|W'|\le4/R$,
\[
|A_t|\le CA+\frac{10\,M_2|\Gamma(b)|}{R}\,\mathbf 1_{\{|L-\lambda_c|\le R\}},
\]
    and the indicator set has $t$-length at most $4R/3$. Since $A\le M_2|\Gamma(b)|$, \eqref{eq:Q} gives 
\[
\int_{J_n}|(A^2)_t|\,dt\le C\|A\|_2^2,
\]
    where $C$ is independent of $\omega,R,n$. Hence, for a constant $C_5$ independent of $\omega,R,n$,
\[
\Bigl|\int_{J_n}A^2\cos2\vartheta\,dt\Bigr|\le\frac{C_5}{\omega}\|A\|_2^2 \qquad\text{and}\qquad \|A\sin\vartheta\|_2^2\ge\frac12\Bigl(1-\frac{C_5}{\omega}\Bigr)\|A\|_2^2.
\]

    \textit{Conclusion.} Using $(p+r)^2\ge p^2/2-r^2$ and \eqref{eq:err},
\[
\|\psi_t\|_2^2\ge\frac{\omega^2}{4}\Bigl(1-\frac{C_5}{\omega}\Bigr)\|A\|_2^2 - 2\|\phi\|_\infty^2\|A\|_2^2-2M_2^2C_\omega
\]
    and
\[
\|\psi\|_2^2\le2\|A\|_2^2+2M_2^2C_\omega.
\]
    Substituting into \eqref{eq:basic},
\[
\int_{S_{n1}}\biggl[\Bigl(\frac{d\varphi_\nu}{ds}\Bigr)^2-\varphi_\nu^2\biggr]ds \ge \|A\|_2^2\Bigl[\frac{\omega^2}{4\sigma_2}\Bigl(1-\frac{C_5}{\omega}\Bigr)-\frac{2\|\phi\|_\infty^2}{\sigma_2}-2\sigma_2\Bigr]-2M_2^2(\sigma_2^{-1}+\sigma_2)C_\omega.
\]
    Since $\phi$, $M_1$, $M_2$, $\sigma_2$, and $C_5$ do not depend on $\omega$, we may choose $\omega=\omega_0$ so that the bracket is at least $1$. Let
\[
R_0 = \max\Bigl\{1,\ \frac{20\,M_2^2(\sigma_2^{-1}+\sigma_2)C_{\omega_0}}{M_1^2|\Gamma(b)|^2}\Bigr\}.
\]
    For $R\ge R_0$, \eqref{eq:Q} yields
\[
\int_{S_{n1}}\biggl[\Bigl(\frac{d\varphi_\nu}{ds}\Bigr)^2-\varphi_\nu^2\biggr]ds \ge \frac25M_1^2|\Gamma(b)|^2R-\frac1{10}M_1^2|\Gamma(b)|^2R.
\]
    Thus, the claim holds with $c_2=3M_1^2|\Gamma(b)|^2/(20\pi)$.
\end{proof}

\begin{lemma}[The remaining sheets]\label{lem:others}
    Let $\omega=\omega_0$, $R\ge R_0$, and let $\lambda_c$ be an interior center. For all sufficiently large $n$ (depending on $R$),
\[
\sum_{j=2}^4\int_{S_{nj}}\varphi_\nu^2\,ds\le\pi c_2R.
\]
    Consequently,
\[
\frac1{2\pi}\sum_{j=2}^4\int_{S_{nj}}\biggl[\Bigl(\frac{d\varphi_\nu}{ds}\Bigr)^2-\varphi_\nu^2\biggr]ds \ge -\frac{c_2}{2}R.
\]
\end{lemma}

\begin{proof}
    The second claim follows from the first, since the derivative terms are nonnegative.

    \textit{The sheet $S_{n2}$.} Note that the reflection $x\mapsto-x$ acts on $g$ through the isometry $\operatorname{diag}((-1)^j)$, so $ds=\sigma(t)\,dt$ on $S_{n2}$ as well. By \eqref{eq:phiS} and \eqref{eq:alt}, $|\varphi_\nu(-x)|\le M_2\bar r(L)$ for $x\in S_{n1}$. Therefore,
\[
\int_{S_{n2}}\varphi_\nu^2\,ds \le \sigma_2M_2^2\|\bar r\|_2^2\le \sigma_2 M_2^2 C_{\omega_0} e^{-2(\lambda_c-R)}\le \sigma_2 M_2^2 C_{\omega_0} e^{-2L(t_n)} \to 0.
\]

    \textit{The sheets $S_{n3}\cup S_{n4}$.} Let $x\in S_{n3}\cup S_{n4}$ and $z=|x|^{-1}\in S_{n1}$, so $z\le1-b_n$. Replacing $\nu$ by $\operatorname{diag}((-1)^j)\nu$ if $x<0$ does not change $\|\nu\|$ or its support, so we may take $x>0$. By \eqref{eq:recip}, $\varphi_\nu(x)=\langle\nu,u\rangle/\sqrt{k^*(z^2)}$ with $u=(v_jz^{n-j})_{j=0}^n$. Since $\nu$ is supported in $\{j\le n/2\}$, the Cauchy--Schwarz inequality gives
\[
|\langle\nu,u\rangle|\le\|\nu\|\,z^{n/2}\Bigl(\sum_{j\le n}v_j^2\Bigr)^{1/2}\le C\|\nu\|(n+1)^{|\tau|+1/2}e^{-d_n/2},
\]
    because $z^{n/2}\le e^{-nb_n/2}=e^{-d_n/2}$. Keeping only the term $j=0$ in $k^*$ gives 
\[
k^*(z^2)\ge v_n^2\ge C_1^2(n+1)^{-2|\tau|}.
\]
    Hence, 
\[
|\varphi_\nu(x)|\le C\|\nu\|(n+1)^{2|\tau|+1/2}e^{-d_n/2}.
\]
    The arclength of $S_{n3}\cup S_{n4}$ is $O(|J_n|)=O(\log n)$ by Section \ref{sec:chaos}, so by \eqref{eq:dn} and \eqref{eq:etanorm},
\[
\int_{S_{n3}\cup S_{n4}}\varphi_\nu^2\,ds \le  C\|\nu\|^2n^{4|\tau|+2}e^{-d_n}\le CC_4Rn^{4|\tau|+2}e^{-d_n}\to0.
\]
    For fixed $R$, both bounds tend to $0$.
\end{proof}

\subsection{Proof of the lower bound} \label{subsec:windows}
    Fix $\omega=\omega_0$ and $R=R_0$ as in Lemma \ref{lem:margin}, and let $c_2$ and $C_4$ be the corresponding constants. Define the centers
\[
\lambda_j=L(t_n)+1+R+(2R+2)(j-1),\qquad 1\le j\le m,
\]
    where $m$ is the largest integer with $\lambda_m+R\le L(u_n)-1$. Each $\lambda_j$ is an interior center, and by \eqref{eq:Ldiff}, for large $n$,
\begin{equation}\label{eq:numwindows}
m\ge\frac{L(u_n)-L(t_n)}{2R+2}-1\ge\frac{\log n}{8R+8}.
\end{equation}

    Let $\nu^{(j)}$ be the vector constructed above with center $\lambda_j$. Since 
\[
e^{\lambda_{j+1}-R}=e^{\lambda_j+R+2}>e^{\lambda_j+R},
\]
    the supports of the $\nu^{(j)}$ are pairwise disjoint, so these vectors are pairwise orthogonal. Let $\mathcal S\in\{\mathcal S_n,\,S_{n1}\cup S_{n3}\}$. By \eqref{eq:quad}, Lemma \ref{lem:margin}, and Lemma \ref{lem:others} (the terms on the sheets $S_{nj}\subset\mathcal S$ with $j\ge2$ are bounded below by $-\frac1{2\pi}\int\varphi_\nu^2\,ds$), for large $n$,
\[
\langle\nu^{(j)},T_{\mathcal S}\nu^{(j)}\rangle \ge  c_2R-\frac{c_2}{2}R=\frac{c_2}{2}R\ge\frac{c_2}{2C_4}\|\nu^{(j)}\|^2=:c_3\|\nu^{(j)}\|^2,
\]
    where \eqref{eq:etanorm} was used in the last inequality. In particular, $\nu^{(j)}\ne0$. Corollary \ref{cor:test} and \eqref{eq:numwindows} give
\[
\Var N_n(\mathcal S)\ge2c_3^2m\ge c\log n,\qquad \mathcal S\in\{\mathcal S_n,\ S_{n1}\cup S_{n3}\}.
\]

    Finally, the map $\xi_j\mapsto(-1)^j\xi_j$ preserves the Gaussian law and sends $P_n(x)$ to $P_n(-x)$. Hence, $(N_n(S_{n2}),N_n(S_{n4}))\overset d=(N_n(S_{n1}),N_n(S_{n3}))$, and so
\[
\Var N_n(S_{n2}\cup S_{n4})=\Var N_n(S_{n1}\cup S_{n3})\ge c\log n.
\]
    This completes the proof of Theorem \ref{thm:gauss}.

\section{The non-centered case}\label{sec:noncentered}
    In this section, the $\xi_j$ are again general random variables satisfying Condition \ref{cond:xi}, $Q_n=M_n+P_n$ is as in Section \ref{subsec:noncentered}, with $|\mu_j|\le C_\mu$, and $N_{Q_n}(S)$ and $N_{P_n}(S)$ are the numbers of real zeros in $S$ of $Q_n$ and $P_n$. Constants may depend on $C_\mu$ and on the parameters $x_0,\alpha,\beta,\gamma$ of Theorem \ref{thm:noncentered}. 

    The following lemma shows that zeros of $Q_n$ away from $\pm1$ contribute $O(1)$ to the second moment. It holds for all real $\mu_j$ with $|\mu_j|\le C_\mu$, in particular, for $\mu_j\equiv0$, where $Q_n=P_n$.

\begin{lemma}[{\cite{Do21}}]\label{lem:O1}
    Fix $x_1\in(0,1)$ and let $A_{x_1}=\{x\in\R:\,||x|-1|\ge x_1\}$. Then $\E N_{Q_n}^2(A_{x_1}) = O(1)$.
\end{lemma}
    A proof, based on an argument of Ibragimov and Maslova, is given in \cite[Lemma 2.2]{Do21}.

    Let $x_0\in(0,1)$ and set $B_+=[x_0,1/x_0]$, $B_-=[-1/x_0,-x_0]$, and $B=B_+\cup B_-$. Since $1/x_0\ge2-x_0$, we have $\R\setminus B\subset A_{1-x_0}$. Lemma \ref{lem:O1} with $x_1=1-x_0$ therefore gives
\begin{equation}\label{eq:outside-B}
\E N_{Q_n}^2(\R\setminus B)=O(1) \qquad\text{and}\qquad \E N_n^2(\R\setminus B)=O(1).
\end{equation}

\begin{lemma}\label{lem:interval}
    Assume Conditions \ref{cond:xi} and \ref{cond:main}. There are constants $\epsilon,C>0$ such that for all sufficiently large $n$ and all $S\in\{\R,B_+,B_-\}$,
\[
\epsilon\log n\le\Var N_n(S)\le C\log n.
\]
\end{lemma}
\begin{proof}
    For $S=\R$, this is Theorem \ref{thm:main}. Let $S=B_+$. Since $\R_+\setminus B_+\subset\R\setminus B$, \eqref{eq:outside-B} gives $\Var N_n(\R_+\setminus B_+)\le C$, and hence $|\sqrt{\Var N_n(B_+)}-\sqrt{\Var N_n(\R_+)}|\le C$. The claim for $B_+$ then follows from Proposition \ref{prop:onesided}. The case $S=B_-$ is identical.
\end{proof}

\subsection{Proof of Theorem \ref{thm:noncentered}}
    We apply Theorem \ref{thm:DNOR} with $d=1/2$.

    \textit{Case \ref{nc:dominated}.} Applying Theorem \ref{thm:DNOR} with $J=S=\R$ gives 
\[
\Var N_{Q_n}=\Var N_n+O(\log^{1/2}n).
\]
    By Lemma \ref{lem:interval}, for large $n$,
\[
\frac12\epsilon\log n \le \Var N_{Q_n}\le C\log n.
\]
    Since $\Var N_n\ge\epsilon\log n$, Theorem \ref{thm:DNOR} also yields the CLT.

    \textit{Case \ref{nc:mixed}.} Assume that $M_n$ is dominated by $P_n$ on $B_-$. The other case reduces to this one by replacing $v_j$ with $(-1)^jv_j$, which replaces $Q_n(x)$ by $Q_n(-x)$, preserves Conditions \ref{cond:xi} and \ref{cond:main} and the bound on $\mu_j$, and interchanges $B_+$ and $B_-$. Since $x_0$ is fixed, the endpoints of $B_\pm$ stay at distance at least $1-x_0$ from $\pm1$, so $J=S=B_\pm$ are admissible choices in Theorem \ref{thm:DNOR}. Applying it with $S=B_+$ and with $S=B_-$ gives
\[
\Var N_{Q_n}(B_+) = O(\log^{1/2}n),\qquad \Var N_{Q_n}(B_-) = \Var N_n(B_-)+O(\log^{1/2}n).
\]
    By Lemma \ref{lem:interval}, 
\[
\frac12\epsilon\log n\le \Var N_{Q_n}(B_-) \le C\log n
\]
    for large $n$. Since $\Var N_n(B_-)\ge\epsilon\log n$, Theorem \ref{thm:DNOR} gives, as $n\to \infty$,
\[
\frac{N_{Q_n}(B_-) - \E N_{Q_n}(B_-)}{\sqrt{\Var N_{Q_n}(B_-)}} \xrightarrow{d}\mathcal N(0,1).
\]
    Let $X=N_{Q_n}(B_+)+N_{Q_n}(\R\setminus B)$, so that $N_{Q_n}=N_{Q_n}(B_-)+X$. By \eqref{eq:outside-B},
\[
\Var X\le2\Var N_{Q_n}(B_+)+2\E N_{Q_n}^2(\R\setminus B)=O(\log^{1/2}n).
\]
    The triangle inequality for the standard deviation gives 
\[
|\sqrt{\Var N_{Q_n}}-\sqrt{\Var N_{Q_n}(B_-)}| = O(\log^{1/4}n).
\]
    Hence, for large $n$, 
\[
c\log n\le\Var N_{Q_n}\le C\log n,
\]
    and as $n\to \infty$,
\[
\frac{\Var N_{Q_n}(B_-)}{\Var N_{Q_n}}\to 1.
\]
    For the CLT, write
\[
\frac{N_{Q_n}-\E N_{Q_n}}{\sqrt{\Var N_{Q_n}}} = \sqrt{\frac{\Var N_{Q_n}(B_-)}{\Var N_{Q_n}}}\ \frac{N_{Q_n}(B_-)-\E N_{Q_n}(B_-)}{\sqrt{\Var N_{Q_n}(B_-)}} + \frac{X-\E X}{\sqrt{\Var N_{Q_n}}}.
\]
    The last term has variance $O(\log^{-1/2}n)$ and therefore tends to $0$ in probability. Slutsky's lemma gives the CLT for $N_{Q_n}$. \qed

\subsection{Verification of Examples \ref{ex:periodic} and \ref{ex:alternating}}
    In both examples, $|v_j|=1$, so $\tau=0$ and Condition \ref{cond:main} holds with $C_1=C_2=1$, and $\mu_j$ is bounded.

    Write $y=1-|x|+1/n$. We will use the following elementary bound. For $|x|\le1$, $l\in\{0,1,2\}$, and $m\ge n-4$,
\begin{equation}\label{eq:elem}
n^l|x|^m\le Cy^{-l}.
\end{equation}
    Indeed, $|x|^m\le e^{-(n-4)(1-|x|)}\le Ce^{-ny}$ and $e^{-s}\le C_ls^{-l}$ for $s=ny\ge1$.

\begin{lemma}\label{lem:variance-profile}
    Let $|v_j|=1$ for $0\le j\le n$. For $l=0,1,2$ and $|x|\le 1$,
\[
\Var P_n^{(l)}(x)=\Var(P_n^*)^{(l)}(x) \asymp y^{-2l-1}.
\]
\end{lemma}
\begin{proof}
    Note that 
\[
\Var P_n^{(l)}(x)=\Var(P_n^*)^{(l)}(x)=\sum_{j=l}^n\bigl((j)_l\bigr)^2x^{2(j-l)}.
\]
    We distinguish three cases.

    If $|x|\le 3/4$, then $y\asymp 1$ and 
\[
1\le(l!)^2\le \Var P_n^{(l)}(x) \le\sum_{j\ge l}(1+j)^{2l}(3/4)^{2(j-l)}\le C.
\]

    If $3/4 \le |x| \le 1-1/n$, then $1-|x|\asymp y$. Since $\bigl((j)_l\bigr)^2\asymp(1+j)^{2l}$ for $j\ge l$, $|x|^{-2l}\asymp1$, and the terms with $j<l$ contribute $O(1)$, we have
\[
\Var P_n^{(l)}(x)\asymp1+\sum_{j=0}^n(1+j)^{2l}x^{2j}.
\]
    By Lemma \ref{lem:sum}, 
\[
\sum_{j=0}^n(1+j)^{2l}x^{2j} \le C(1-x^2)^{-2l-1}\asymp y^{-2l-1}.
\]
    Restricting to $j\le(1-|x|)^{-1}\le n$, where $x^{2j}\ge c$, shows that 
\[
\sum_{j=0}^n(1+j)^{2l}x^{2j} \ge c(1-|x|)^{-2l-1}\asymp y^{-2l-1}.
\]

    If $1-1/n<|x|\le 1$, then $x^{2(j-l)} \asymp 1$ for $l\le j\le n$, and 
\[
\Var P_n^{(l)}(x) \asymp \sum_{j=l}^n(1+j)^{2l} \asymp n^{2l+1} \asymp y^{-2l-1}.
\]
\end{proof}

\begin{proof}[Verification of Example \ref{ex:periodic}]
    Let $K=\lfloor(n+1)/4\rfloor$, so $4K\in[n-2,n+1]$. Since $1+x-x^2-x^3=(1+x)(1-x^2)$, summing over blocks of four consecutive indices gives
\[
M_n(x) = \mu\Bigl[\frac{1+x}{1+x^2}\bigl(1-x^{4K}\bigr) + \sum_{j=4K}^nv_jx^j\Bigr],
\]
    and the last sum has at most three terms. The function $(1+x)/(1+x^2)$ has bounded derivatives on $[-1,1]$, while for $4K\le j\le n$,
\[
\frac{d^l}{dx^l}x^{4K} = O(n^l|x|^{4K-l}) \qquad\text{and}\qquad \frac{d^l}{dx^l}x^j=O(n^l|x|^{j-l}).
\]
    Since $4K-l\ge n-4$ and $j-l\ge n-4$, \eqref{eq:elem} gives
\[
|M_n^{(l)}(x)|\le C|\mu|\,y^{-l},\qquad|x|\le1,\ l=0,1,2.
\]
    The same holds for $M_n^*$. Indeed, the sum $\sum_{j=0}^n v_{n-j}x^j$ over four consecutive indices equals $\pm(1\pm x)(1-x^2)$ for each residue of $n$ modulo $4$, so $M_n^*$ has the same structure. By Lemma \ref{lem:variance-profile}, 
\[
\sqrt{\Var P_n^{(l)}(x)}\ge cy^{-l-1/2},
\]
    and therefore
\[
|M_n^{(l)}(x)| \le C|\mu|\, y^{1/2} \sqrt{\Var P_n^{(l)}(x)},\qquad |x|\le 1,\ l=0,1,2,
\]
    and likewise for $(M_n^*)^{(l)}$ and $(P_n^*)^{(l)}$. Thus, $M_n$ is dominated by $P_n$ on $\R$ up to order $2$ with factor function $\Phi_{\beta,1/2}$ and $\beta=C|\mu|$.
\end{proof}

\begin{proof}[Verification of Example \ref{ex:alternating}]
    Here 
\[
M_n(x)=M_n^*(x)=\mu \sum_{j=0}^n x^j \qquad\text{and}\qquad \Var P_n(x)= \Var P_n^*(x) = \sum_{j=0}^n x^{2j},
\]
    so it suffices to treat $|x|\le 1$.

    For $x\in [-1,0)$, $M_n(x)=\mu(1-x^{n+1})/(1-x)$ with $1-x>1$. Differentiating, we get 
\[
|M_n(x)|\le2|\mu|,\qquad |M_n'(x)|\le|\mu|(2+(n+1)|x|^n),
\]
    and 
\[
|M_n''(x)|\le|\mu|(4+2(n+1)|x|^n+n(n+1)|x|^{n-1}).
\]
    By \eqref{eq:elem}, 
\[
|M_n^{(l)}(x)|\le C|\mu|y^{-l},\qquad x\in [-1,0),\, l=0,1,2.
\]
    As in the previous example, Lemma \ref{lem:variance-profile} shows that $M_n$ is dominated by $P_n$ on $[-1,0)$ up to order $2$ with factor function $\Phi_{\beta,1/2}$, where $\beta=C|\mu|$.

    Assume now that $x\in (0,1]$. Since $x^{2j}\le x^j$ for $0\le j\le n$, Lemma \ref{lem:variance-profile} gives, for some constant $c>0$,
\[
|M_n(x)|\ge |\mu| \Var P_n(x) \ge c|\mu|\, y^{-1/2} \sqrt{\Var P_n(x)}.
\]
    Fix $\alpha>0$. The function $y\mapsto y|\log y|$ increases on $(0,e^{-1})$ and tends to $0$ as $y\to0$, so there is $y_0\in(0,e^{-1}]$, depending only on $\alpha$ and $c|\mu|$, such that 
\[
c|\mu| (y|\log y|)^{-1/2} \ge \alpha,\qquad 0<y\le y_0.
\]
    Let $x_0=1-y_0/2$. Then for all $n\ge 2/y_0$ and all $x\in [x_0,1]$, we have $y\le y_0$ and
\[
|M_n(x)|\ge \alpha\, |\log y|^{1/2} \sqrt{\Var P_n(x)},
\]
    which is the claimed domination with factor function $\Phi_\alpha$.
\end{proof}


\begin{thebibliography}{99}
\bibitem{AN26}
    A. Aguirre and H. H. Nguyen,
    Central limit theorem for real zeros of random Weyl polynomials with general coefficients,
    \href{https://arxiv.org/abs/2609.06754}{arXiv:2609.06754}, 2026.

\bibitem{ANW26}
    A. Aguirre, H. H. Nguyen, and J. Wang,
    Real roots of random Weyl polynomials with general coefficients: expectation and variance,
    \textit{Electron. J. Probab.} \textbf{31} (2026), 52 pp. \MR{5113821} 

\bibitem{AL21}
    M. Ancona and T. Letendre,
    Roots of Kostlan polynomials: moments, strong law of large numbers and central limit theorem,
    \textit{Ann. H. Lebesgue} \textbf{4} (2021), 1659--1703. \MR{4353974}

\bibitem{AL21b}
    M. Ancona and T. Letendre,
    Zeros of smooth stationary Gaussian processes,
    \textit{Electron. J. Probab.} \textbf{26} (2021), Paper No. 68, 81 pp. \MR{4262341}

\bibitem{ABF23}
    E. Assaf, J. Buckley, and N. D. Feldheim,
    An asymptotic formula for the variance of the number of zeroes of a stationary Gaussian process,
    \textit{Probab. Theory Related Fields} \textbf{187} (2023), no. 3--4, 999--1036. \MR{4664589}

\bibitem{ADL16}
    J.-M. Aza\"is, F. Dalmao, and J. Le\'on,
    CLT for the zeros of classical random trigonometric polynomials,
    \textit{Ann. Inst. Henri Poincar\'e Probab. Stat.} \textbf{52} (2016), no. 2, 804--820. \MR{3498010}

\bibitem{AL13}
    J.-M. Aza\"is and J. Le\'on,
    CLT for crossings of random trigonometric polynomials,
    \textit{Electron. J. Probab.} \textbf{18} (2013), no. 68, 17 pp. \MR{3084654}

\bibitem{AW09} 
    J.-M. Aza\"is and M. Wschebor, 
    \textit{Level sets and extrema of random processes and fields}, 
    John Wiley \& Sons, Hoboken, NJ, 2009. \MR{2478201}

\bibitem{BCP19}
    V. Bally, L. Caramellino, and G. Poly,
    Non universality for the variance of the number of real roots of random trigonometric polynomials,
    \textit{Probab. Theory Related Fields} \textbf{174} (2019), no. 3--4, 887--927. \MR{3980307}

\bibitem{BD97}
    P. Bleher and X. Di,
    Correlations between zeros of a random polynomial,
    \textit{J. Stat. Phys.} \textbf{88} (1997), no. 1--2, 269--305. \MR{1468385}

\bibitem{BP32} 
    A. Bloch and G. P\'olya, 
    On the roots of certain algebraic equations, 
    \textit{Proc. London Math. Soc. (2)} \textbf{33} (1932), 102--114. \MR{1576817}

\bibitem{Dal15}
    F. Dalmao,
    Asymptotic variance and CLT for the number of zeros of Kostlan Shub Smale random polynomials,
    \textit{C. R. Math. Acad. Sci. Paris} \textbf{353} (2015), no. 12, 1141--1145. \MR{3427922}

\bibitem{Do21} 
    Y. Do, 
    Real roots of random polynomials with coefficients of polynomial growth: a comparison principle and applications, 
    \textit{Electron. J. Probab.} \textbf{26} (2021), Paper No. 144, 45 pp. \MR{4346676}

\bibitem{DNN22}
    Y. Do, H. H. Nguyen, and O. Nguyen,
    Random trigonometric polynomials: universality and non-universality of the variance for the number of real roots,
    \textit{Ann. Inst. Henri Poincar\'e Probab. Stat.} \textbf{58} (2022), no. 3, 1460--1504. \MR{4452640}

\bibitem{DNNP24}
    Y. Do, H. H. Nguyen, O. Nguyen, and I. E. Pritsker,
    Central limit theorem for the number of real roots of random orthogonal polynomials,
    \textit{Ann. Inst. Henri Poincar\'e Probab. Stat.} \textbf{60} (2024), 2211--2240. \MR{4780518}

\bibitem{DHNV15} 	    
    Y. Do, H. H. Nguyen, and V. Vu,
    Real roots of random polynomials: expectation and repulsion,
    \textit{Proc. Lond. Math. Soc. (3)} \textbf{111} (2015), no. 6, 1231--1260. \MR{3447793}

\bibitem{DN25} 
    Y. Do and N. D. V. Nguyen, 
    Real roots of random polynomials: asymptotics of the variance, 
    \textit{Electron. J. Probab.} \textbf{30} (2025), article no. 112, 37 pp. \MR{4939546}

\bibitem{DNOR26} 
    Y. Do, N. D. V. Nguyen, and S. O'Rourke, 
    Real roots of non-centered random polynomials, \href{https://arxiv.org/abs/2605.26402}{arXiv:2605.26402}, 2026.

\bibitem{DNV18} 	    
    Y. Do, O. Nguyen, and V. Vu, 
    Roots of random polynomials with coefficients of polynomial growth,
    \textit{Ann. Probab.} \textbf{46} (2018), no. 5, 2407--2494. \MR{3846831}

\bibitem{DV20}
    Y. Do and V. Vu,
    Central limit theorems for the real zeros of Weyl polynomials,
    \textit{Amer. J. Math.} \textbf{142} (2020), no. 5, 1327--1369. \MR{4150647}

\bibitem{EK95} 
    A. Edelman and E. Kostlan,
    How many zeros of a random polynomial are real?, 
    \textit{Bull. Amer. Math. Soc. (N.S.)} \textbf{32} (1995), 1--37. \MR{1290398}

\bibitem{EO56} 
    P. Erd\H{o}s and A. C. Offord, 
    On the number of real roots of a random algebraic equation, 
    \textit{Proc. London Math. Soc. (3)} \textbf{6} (1956), 139--160. \MR{0073870}

\bibitem{Gas23}
    L. Gass,
    Variance of the number of zeros of dependent Gaussian trigonometric polynomials,
    \textit{Proc. Amer. Math. Soc.} \textbf{151} (2023), 2225--2239. \MR{4556213}

\bibitem{GW11}
    A. Granville and I. Wigman,
    The distribution of the zeros of random trigonometric polynomials,
    \textit{Amer. J. Math.} \textbf{133} (2011), no. 2, 295--357. \MR{2797349}

\bibitem{IM71P1} 
    I. A. Ibragimov and N. B. Maslova, 
    The mean number of real zeros of random polynomials. I. Coefficients with zero means, 
    \textit{Teor. Verojatnost. i Primenen.} \textbf{16} (1971), 229--248. \MR{0286157}

\bibitem{IM71II}
    I. A. Ibragimov and N. B. Maslova,
    The mean number of real zeros of random polynomials. II. Coefficients with a nonzero mean,
    \textit{Teor. Verojatnost. i Primenen.} \textbf{16} (1971), 495--503. \MR{0288824} 

\bibitem{Kac43} 
	M. Kac, 
	On the average number of real roots of a random algebraic equation,
	\textit{Bull. Amer. Math. Soc.} \textbf{49} (1943), 314--320. \MR{0007812}

\bibitem{Kac49} 
	M. Kac, 
	On the average number of real roots of a random algebraic equation II,
	\textit{Proc. London Math. Soc. (2)} \textbf{50} (1949), 390--408. \MR{0030713}

\bibitem{LO43} 
	J. E. Littlewood and A. C. Offord, 
	On the number of real roots of a random algebraic equation III,
	\textit{Rec. Math. [Mat. Sbornik] N.S.} \textbf{12(54)} (1943), 277--286. \MR{0009656}

\bibitem{LO45} 
	J. E. Littlewood and A. C. Offord,
	On the distribution of the zeros and $a$-values of a random integral function I,
	\textit{J. London Math. Soc.} \textbf{20} (1945), 130--136. \MR{0019123}

\bibitem{LO48} 
	J. E. Littlewood and A. C. Offord,
	On the distribution of zeros and $a$-values of a random integral function II,
	\textit{Ann. of Math. (2)}  \textbf{49} (1948), 885--952; errata \textbf{50} (1949), 990--991. \MR{0029981}

\bibitem{LP21}
    D. S. Lubinsky and I. E. Pritsker,
    Variance of real zeros of random orthogonal polynomials,
    \textit{J. Math. Anal. Appl.} \textbf{498} (2021), Paper No. 124954, 32 pp. \MR{4202193}

\bibitem{LP22}
    D. S. Lubinsky and I. E. Pritsker,
    Variance of real zeros of random orthogonal polynomials for varying and exponential weights,
    \textit{Electron. J. Probab.} \textbf{27} (2022), Paper No. 83, 32 pp. \MR{4444378}

\bibitem{M74I} 
    N. B. Maslova, 
    The distribution of the number of real roots of random polynomials, 
    \textit{Teor. Verojatnost. i Primenen.} \textbf{19} (1974), 488--500. \MR{0368136}

\bibitem{M74} 
    N. B. Maslova, 
    The variance of the number of real roots of random polynomials, 
    \textit{Teor. Verojatnost. i Primenen.} \textbf{19} (1974), 36--51. \MR{0334327}

\bibitem{NNV16}     
    H. H. Nguyen, O. Nguyen, and V. Vu, 
    On the number of real roots of random polynomials,
    \textit{Commun. Contemp. Math.} \textbf{18} (2016), no. 4, 1550052, 17 pp. \MR{3493213}

\bibitem{N24}
    N. D. V. Nguyen,
    The number of real zeros of elliptic polynomials,
    \textit{Electron. J. Probab.} \textbf{29} (2024), Paper No. 81, 49 pp. \MR{4757541}

\bibitem{NV21} 
    O. Nguyen and V. Vu,
    Random polynomials: central limit theorems for the real roots, 
    \textit{Duke Math. J.} \textbf{170} (2021), no.~17, 3745--3813. \MR{4340724}

\bibitem{NV22} 
    O. Nguyen and V. Vu,
    Roots of random functions: a framework for local universality, 
    \textit{Amer. J. Math.} \textbf{144} (2022), 1--74. \MR{4367414}

\bibitem{R87} 
    W. Rudin, 
    \textit{Real and complex analysis}, 3rd ed., McGraw-Hill, New York, 1987. \MR{0924157}

\bibitem{SM08}
G. Schehr and S. N. Majumdar,
Real roots of random polynomials and zero crossing properties of diffusion equation,
\textit{J. Stat. Phys.} \textbf{132} (2008), no. 2, 235--273. \MR{2415102}

\bibitem{St69} 	   
    D. C. Stevens, 
    The average number of real zeros of a random polynomial,
    \textit{Comm. Pure Appl. Math.} \textbf{22} (1969), 457--477. \MR{0251003}

\bibitem{TV15}
    T. Tao and V. Vu,
    Local universality of zeroes of random polynomials,
    \textit{Int. Math. Res. Not. IMRN} \textbf{2015}, no. 13, 5053--5139. \MR{3439098}
\end{thebibliography}
\end{document}